\documentclass[a4paper,reqno,12pt]{amsart}
\usepackage{amsmath,amstext,amssymb,amsopn,amsthm,mathrsfs}
\usepackage{xcolor}
\numberwithin{equation}{section}
\allowdisplaybreaks

\newtheorem{theorem}{Theorem}[section]
\newtheorem{proposition}[theorem]{Proposition}
\newtheorem{lemma}[theorem]{Lemma}
\newtheorem{corollary}[theorem]{Corollary}
\newtheorem{definition}[theorem]{Definition}

\newcommand{\N}{\mathbb{N}}
\newcommand{\NN}{\mathcal{N}}
\newcommand{\C}{\mathbb{C}}
\newcommand{\Om}{\Omega}
\newcommand{\Omp}{\Omega_+}
\newcommand{\R}{\mathbb{R}^d}
\newcommand{\Rp}{\mathbb{R}^d_+}
\newcommand{\RR}{\mathbb{R}^2}
\newcommand{\Zd}{\mathbb{Z}^d_2}
\newcommand{\EE}{\mathcal{E}}
\newcommand{\EEeta}{\mathcal{E}^\eta}
\newcommand{\A}{\mathcal{A}}
\newcommand{\D}{{\rm Dom}}

\def\a{\alpha}
\def\b{\beta}
\def\la{\lambda}
\def\v{\varphi}
\def\s{\sigma}
\def\om{\omega}
\def\ve{\varepsilon}

\begin{document}
\title[Finite reflection groups and symmetric extensions of Laplacian]{ Finite reflection groups and symmetric extensions of Laplacian }

\author[K. Stempak]{Krzysztof Stempak}
\address{Wydzia\l{} Matematyki \\
         Politechnika Wroc\l{}awska\\ 
         Wyb{.} Wyspia\'nskiego 27\\ 
         50-370 Wroc\l{}aw, Poland}        
\email{Krzysztof.Stempak@pwr.edu.pl}

\date{}

\begin{abstract} Let $W$ be a finite reflection group associated with a root system $R$ in $\mathbb R^d$. Let $C_+$ denote a positive Weyl chamber. 
Consider an open subset $\Omega$ of $\mathbb R^d$, symmetric with respect to reflections from $W$. Let $\Omega_+=\Omega\cap C_+$ be the positive 
part of $\Omega$. We define a family $\{-\Delta_{\eta}^+\}$ of self-adjoint extensions of the Laplacian $-\Delta_{\Omp}$, labeled by homomorphisms $\eta\colon W\to \{1,-1\}$. In the construction of these $\eta$-Laplacians  $\eta$-symmetrization of functions on $\Omega$ is involved. The Neumann Laplacian 
$-\Delta_{N,\Omega_+}$ is included and corresponds to $\eta\equiv1$. If $H^{1}(\Om)=H^{1}_0(\Om)$, then the Dirichlet Laplacian $-\Delta_{D,\Omega_+}$ is 
either included and corresponds to $\eta={\rm sgn}$; otherwise the Dirichlet Laplacian is considered separately. Applying  the spectral functional calculus we consider the pairs of operators $\Psi(-\Delta_{N,\Omega})$ and $\Psi(-\Delta_{\eta}^+)$, or  $\Psi(-\Delta_{D,\Omega})$ and $\Psi(-\Delta_{D,\Omega_+})$, 
where  $\Psi$ is a Borel function on $[0,\infty)$. We prove relations between the integral kernels for the operators in these pairs, which are given in 
terms of symmetries governed by $W$. 
\end{abstract}

\subjclass[2010]{Primary 35K08; Secondary 47B25.}

\keywords{ Neumann Laplacian, Dirichlet Laplacian, self-adjoint operator, finite reflection group, Sobolev space, sesquilinear form, functional calculus, 
heat kernel.}

\maketitle

\section{Introduction} \label{sec:intro}

Let $\Om$ be  a nonempty open subset of $\mathbb R^d$, $d\ge1$, and let $\Delta=\sum_1^d \partial_j^2$ denote the Laplacian. If not otherwise stated, 
$-\Delta_\Om$ will mean the differential operator $f\mapsto -\Delta f$ with dense in $L^2(\Om)$ domain $C^\infty_c(\Om)$. Clearly  $-\Delta_\Om$ is symmetric,
$$
\langle(-\Delta_\Om)f,g\rangle_{L^2(\Om)}=\langle f,(-\Delta_\Om)g\rangle_{L^2(\Om)}, \qquad f,g\in {\rm Dom}(-\Delta_\Om)=C^\infty_c(\Om),
$$
and non-negative, $\langle (-\Delta_\Om)f,f\rangle_{L^2(\Om)}\ge0$ for $f\in {\rm Dom}(-\Delta_\Om)$. 

Let $\mathfrak{t}_\Om$ be the sesquilinear form defined on the Sobolev space $H^{1}(\Omega)$ as its domain by
$$
\mathfrak{t}_\Om[f,g]=\int_{\Om}(\nabla f)(x)\cdot\overline{(\nabla g)(x)}\,dx=\int_{\Om}\sum_{j=1}^d \partial_jf(x)\,\overline{\partial_jg(x)}\,dx. 
$$
The \textit{Neumann Laplacian} on $\Om$, denoted $-\Delta_{N,\,\Om}$, is defined as the operator on $L^2(\Omega)$ associated with the  
form $\mathfrak{t}_{N,\,\Omega}:=\mathfrak{t}_\Om$; in particular, ${\rm Dom}(-\Delta_{N,\,\Om})\subset{\rm Dom}(\mathfrak{t}_{N,\,\Omega}):= H^{1}(\Omega)$. 
On the other hand, the \textit{Dirichlet Laplacian} on $\Om$, denoted $-\Delta_{D,\,\Om}$, is defined as the operator on $L^2(\Omega)$ associated with the 
form $\mathfrak{t}_{D,\,\Omega}$, which is the restriction of $\mathfrak{t}_{\Omega}$ to $H^{1}_0(\Om)$; in particular, 
${\rm Dom}(-\Delta_{D,\,\Om})\subset{\rm Dom}(\mathfrak{t}_{D,\,\Omega}):= H^{1}_0(\Omega)$.
Since the forms $\mathfrak{t}_{N,\Omega}$ and $\mathfrak{t}_{D,\Omega}$ are Hermitian, closed and non-negative, the associated operators are self-adjoint 
and non-negative. See \cite[Chapter 10 and Section 3 of Chapter 12]{Sch} and also \cite[pp. 263 and 265]{RS}, where comments on the definitions of 
$-\Delta_{N/D,\,\Om}$ are gathered. Each of the operators $-\Delta_{N/D,\Om}$ is indeed an extension of $-\Delta_{\Om}$. This follows from the definitions 
in terms of forms, with an application of Gauss' formula for functions from Sobolev classes (that can be found, for instance, in \cite[Appendix D]{Sch}; 
see comments at the end of Section \ref{sec:ext}). We also mention that $-\Delta_{D,\Om}$ coincides with the Friedrichs extension of 
$\overline{-\Delta_{\Om}}$, the closure of $-\Delta_{\Om}$. See \cite[Section 10.6.1]{Sch}.

In the setting of a general open set $\Om$ it is known (see, for instance, \cite[Section 10.6.1]{Sch}) that
\begin{equation*}
{\rm Dom}(-\Delta_{D,\,\Om,})=H^\Delta(\Om)\cap H^{1}_0(\Om),
\end{equation*}
where $H^\Delta(\Om)=\{f\in L^2(\Om)\colon \Delta f\in L^2(\Om)\}$ and for $f\in L^2(\Om)\subset C^\infty_c(\Om)'$, $\Delta f$ is understood in the 
distributional sense, and
\begin{equation*}
-\Delta_{D,\,\Om}f=-\Delta f\,\,\,\,{\rm for}\,\,\, f\in {\rm Dom}(-\Delta_{D,\,\Om}).
\end{equation*}
Note that $H^{2}(\Om)\subset H^\Delta(\Om)$ but in general the inclusion may be proper. Contrary to the case of the Dirichlet Laplacian much less is known 
about the explicit description of ${\rm Dom}(-\Delta_{N,\,\Om})$, the domain of the Neumann Laplacian, in the setting of general $\Om\subset \R$.

If $\Om$ is an open bounded subset in $\R$, $d\ge2$, with boundary $\partial\Om$ of class $C^2$, or an open bounded subset of $\mathbb{R}$, then there are 
much finer results concerning properties of $-\Delta_{D,\,\Omega}$ and $-\Delta_{N,\,\Omega}$. In particular, in this case the Dirichlet Laplacian refers to vanishing boundary values at $\partial\Om$ and the Neumann Laplacian refers to vanishing directional normal derivatives at $\partial\Om$,  see for instance 
\cite[Theorems 10.19 and 10.20]{Sch}.

In the special case $\Om=\R$, one has $H^{1}(\R)= H^{1}_0(\R)$ and it follows that $-\Delta_{N,\,\R}=-\Delta_{D,\,\R}$. More precisely,
$$
{\rm Dom}(-\Delta_{N/D,\,\R})=H^{2}(\R) \,\,{\rm and}\quad -\Delta_{N/D,\,\R}f=-\Delta f\,\,\,{\rm for}\,\,\, f\in H^{2}(\R).
$$

We refer the reader to \cite{AF} for a comprehensive treatment of the theory of Sobolev spaces.

By the spectral theorem, for a Borel function $\Psi$ on $[0,\infty)$, we consider the operators $\Psi(-\Delta_{N/D,\,\Om})$ (recall that $-\Delta_{N/D,\,\Om}$ 
are non-negative and hence their spectra are contained in $[0,\infty)$). In particular, we associate with the Dirichlet Laplacian $-\Delta_{D,\,\Om}$ the 
semigroup  $\{\exp(-t(-\Delta_{D,\,\Om}))\}_{t>0}$ of bounded on $L^2(\Om)$ operators, called the \textit{Dirichlet heat semigroup}.
Each  $\exp(-t(-\Delta_{D,\,\Om}))$, $t>0$, is an integral operator with a kernel $p_t^{D,\,\Om}(x,y)$, that is for every $f\in L^2(\Om)$ there holds
$$
\exp(-t(-\Delta_{D,\,\Om}))f(x)=\int_{\Om}p_t^{D,\,\Om}(x,y)f(y)dy, \qquad x-a.e..
$$
Moreover, as a function on $(0,\infty)\times\Om\times\Om$, $p_t^{D,\,\Om}(x,y)$ is $C^\infty$ and strictly positive. See \cite[Theorem 5.2.1]{D1}.
Then $\{p_t^{D,\,\Om}(x,y)\}_{t>0}$,  is called the \textit{Dirichlet heat kernel on} $\Om$.

Analogously, we consider  $\{\exp(-t(-\Delta_{N,\,\Om}))\}_{t>0}$, the \textit{Neumann heat semigroup} associated with  $-\Delta_{N,\,\Om}$. As before,
each  $\exp(-t(-\Delta_{N,\,\Om}))$, $t>0$, is an integral operator with a kernel $p_t^{N,\,\Om}(x,y)$
which, as a function on $(0,\infty)\times\Om\times\Om$, is $C^\infty$ and strictly positive. 
Then $\{p_t^{N,\,\Om}(x,y)\}_{t>0}$,  is called the \textit{Neumann heat kernel on} $\Om$.

Let  $R$ be a \textit{normalized root system} in $\R$, that is a finite set of unit vectors such that $\s_\a(R)=R$ for every $\a\in R$, 
where $\s_\a$ denotes the orthogonal reflection in $\langle \a\rangle^\bot$, the hyperplane orthogonal to $\a$,
$$
\s_\a(x)=x-2\langle \a,x\rangle\a,  \qquad x\in\R.
$$
Clearly, $R\cap\mathbb{R}\a=\{\a,-\a\}$ for every $\a\in R$. The \textit{finite reflection group} $W=W(R)$ \textit{associated with} $R$  
(the group is indeed finite as simple arguments show) is the subgroup of $O(\mathbb R^d)$ generated by the reflections $\s_\a$, 
$\a\in R$. The set $\R\setminus\bigcup_{\a\in R}\langle \a\rangle^\bot$ splits into an even number (equal to $|W|$) of connected components which 
are open polyhedral cones called the \textit{Weyl chambers}. $W$ acts on the set of Weyl chambers, the action is simply transitive, and hence the 
Weyl chambers are mutually congruent. A choice of $\check\a\in \R$ such that  $\langle \a,\check\a\rangle\neq0$ for every $\a\in R$, gives the partition 
$R=R_+\sqcup(- R_+)$, where $R_+=\{\a\in R\colon \langle \a,\check\a\rangle>0\}$. $R_+$ is then referred to as the \textit{set of positive roots}. The 
partition distinguishes the chamber $C_+=\{x\in\Om\colon \forall\, \a\in R_+ \,\,\,\langle x,\a\rangle>0\}$, which is called the \textit{positive Weyl chamber}. 

For a comprehensive treatment of the general theory of finite reflection groups the reader is kindly referred to \cite{Hu}, \cite{He} and \cite[Chapter 4]{DX}.

Given  $\Om\subset\R$ we say it is $W$-\,symmetric provided $\s_\a(\Om)=\Om$ for $\a\in R$ (which implies that $g(\Om)=\Om$ for $g\in W$). 
Then we distinguish the \textit{positive part of} $\Om$ by setting $\Omp:=\Om\cap C_+$.

In the particular case of $\Om=\R$ and $R_+$ consisting of a single vector $\a$, say $\a=(0,\ldots,0,1)$, so that $\Omp=\{x\colon x_d>0\}$ is the half-space 
$\mathbb R^{d-1}\times (0,\infty)$, the following \textit{reflection principles} for the Neumann and Dirichlet heat kernels on the half-space are well-known 
('+' for N[eumann], '$-$' for D[irichlet])
$$
p_t^{N/D,\Omp}(x,y)=p_t(x,y)\pm p_t(\tilde x,y), \qquad x,y\in \Omp.
$$
Here $\{p_t\}_{t>0}$ denotes the Gauss-Weierstrass kernel on $\R$, $\{p_t^{N/D,\Omp}\}_{t>0}$ denote the Neumann and Dirichlet heat kernels on $\Omp$, 
respectively, and $\tilde x=(x_1,x_2,\ldots, -x_d)$ is the reflection point of $x=(x_1,x_2,\ldots, x_d)$ in the hyperplane $x_d=0$. In \cite{MS} these formulas 
were proved to hold in a general context of an arbitrary open subset of $\R$, symmetric in $\langle\a\rangle^\bot$. Moreover, the case of a finite number of symmetries with respect to mutually orthogonal hyperplanes was also considered and similar formulas with natural modifications  were obtained by a recursion argument from the case of a single reflection. To be precise heat kernels were embedded in a more general context of action of functions on operators by a use 
of spectral functional calculus. In its simplest form, the reflection principles extracted from  \cite[Corollary 4.1]{MS} for the heat kernels 
$\{p_t^{N/D,\,\Om_{+}}\}_{t>0}$, corresponding to the Neumann and Dirichlet Laplacians on $\Om_{+}=\Om\cap \R_{+}$, the positive part of an open subset 
$\Om$ of $\R$ symmetric in $\langle e_j\rangle^\bot$, $j=1,\ldots,d$, where  $\R_+=(0,\infty)^d$, read as follows:
\begin{equation*}
p_t^{N,\Omp}(x,y)=\sum_{\varepsilon\in\{-1,1\}^d}p_t^{N,\Om}(\varepsilon x,y),\qquad x,y\in\Om_{+},
\end{equation*} 
\begin{equation*}
p_t^{D,\Omp}(x,y)=\sum_{\varepsilon\in\{-1,1\}^d}{\rm sgn}(\varepsilon)p_t^{D,\Om}(\varepsilon x,y),\qquad x,y\in\Om_{+}.
\end{equation*} 
Here for $x=(x_1,\ldots,x_d)\in \R_{+}$ and $\varepsilon\in\{-1,1\}^d$ 
we write ${\rm sgn}(\varepsilon)=\prod_{i=1}^d\varepsilon_i$  and  $\varepsilon x=(\varepsilon_1 x_1,\ldots \varepsilon_d x_d)$.

It became then clear that the results obtained for a single reflection or in a slightly more general setting of reflections with respect to mutually orthogonal hyperplanes should find their final and complete destination in the context of symmetries governed by an arbitrary finite reflection group. This was, in fact, the main motivation for conducting the present research. However, a closer insight into the problem revealed that using constant or alternating signs to label summands on the right hand sides of expressions for $p_t^{N/D,\Omp}(x,y)$  meant using one of the two special homomorphisms of $W\simeq\{-1,1\}^d$  into 
$\widehat{\mathbb Z}_2=\{-1,1\}$ with multiplication. This observation motivated us strongly to  introduce and investigate a family of natural self-adjoint extensions of $-\Delta_{\Omp}$, the Laplacian on $\Omp$,  by means of homomorphisms $\eta\in{\rm Hom}(W,\widehat{\mathbb Z}_2)$, and to connect them in terms of functional calculus to the Neumann  Laplacian on $\Om$. The main ingredient of the construction of the $\eta$-Laplacian $-\Delta^+_\eta$  consist in $\eta$-symmetrization of functions on $\Om$. See Section \ref{sec:ext} for definitions. For $\eta=\textbf{1}$ the corresponding $\textbf{1}$-Laplacian is the Neumann Laplacian on $\Omp$. For $\eta={\rm sgn}$, in the case when $H^{1}_0(\Om)=H^{1}(\Om)$, the corresponding ${\rm sgn}$-Laplacian is the Dirichlet Laplacian on $\Omp$. Here and later on ${\rm sgn}$ denotes the function ${\rm sgn}(g):={\rm det}(g)\in\{\pm 1\}$ on $W$. The case when $H^{1}_0(\Om)\neq H^{1}(\Om)$ (for instance, for any bounded $\Om$) is treated separately and then analogous connections between $-\Delta_{D,\,\Omp}$ and $-\Delta_{D,\,\Om}$, the Dirichlet Laplacians on $\Omp$ and $\Om$, respectively, are established. 

We prove that the integral kernels of operators emerging in spectral calculus applied to the $\eta$-Laplacians or Dirichlet Laplacians on $\Omp$ in the general 
setting of an open $W$-\,symmetric $\Om\subset\R$, where $W\subset O(\R)$ is an arbitrary finite reflection group, are related to the corresponding integral 
kernels of Neumann or Dirichlet Laplacians on $\Om$. More specifically, we prove the following.
\begin{theorem} \label{thm:main} Let $\Om$ be an open  $W$-\,symmetric subset of $\R$  with $\Om_+$ as its positive part. Let $\Psi$ be a Borel 
function on $[0,\infty)$ and $\eta\in{\rm Hom}(W,\,\widehat{\mathbb Z}_2)$. Assume that $\Psi(-\Delta_{N,\,\Om})$ is an integral operator with  
kernel $K^\Psi_{-\Delta_{N,\,\Om}}$. Then $\Psi(-\Delta^+_\eta)$ is also an integral operator with kernel $K^\Psi_{-\Delta^+_\eta}$ given by
\begin{equation}\label{kerN}
K^\Psi_{-\Delta^+_\eta}(x,y)= \sum_{g\in W}\eta(g)K^\Psi_{-\Delta_{N,\,\Om}}(gx,y),\qquad x,y\in\Om_+.
\end{equation}
Similarly, if $\Psi(-\Delta_{D,\,\Om})$ is an integral operator with  kernel $K^\Psi_{-\Delta_{D,\,\Om}}$, then $\Psi(-\Delta_{D,\,\Om_+})$ is also an integral operator with kernel $K^\Psi_{-\Delta_{D,\,\Om_+}}$ given by
\begin{equation}\label{kerD}
K^\Psi_{-\Delta_{D,\,\Om_+}}(x,y)=\sum_{g\in  W}{\rm sgn}(g)\,K^\Psi_{-\Delta_{D,\,\Om}}(gx,y),\qquad x,y\in\Om_+.
\end{equation}
\end{theorem}

As a direct corollary of Theorem \ref{thm:main} we obtain the following result for heat kernels.  
\begin{corollary} \label{cor:main}
Let $\Om$ be an open  $W$-\,symmetric subset of $\R$ and let  $p_t^{\eta,\,\Om_+}$ and $p_t^{D,\,\Om_+}$, and $p_t^{N,\,\Om}$ and $p_t^{D,\,\Om}$, 
denote the $\eta$ and the Dirichlet heat kernels on $\Om_+$ and the Neumann and the Dirichlet heat kernels on $\Om$, respectively. Then
\begin{equation*}
p_t^{\eta,\,\Om_+}(x,y)= \sum_{g\in W}\eta(g)p_t^{N,\,\Om}(gx,y) ,\qquad x,y\in\Om_+,\quad t>0,
\end{equation*}
and
\begin{equation*}
p_t^{D,\,\Om_+}(x,y)= \sum_{g\in W}{\rm sgn}(g)\,p_t^{D,\,\Om}(gx,y),\qquad x,y\in\Om_+,\quad t>0.
\end{equation*}
\end{corollary}
We point out that analogous formulas hold for other important kernels like resolvent kernels, Riesz potential kernels and so on.

The essential novelty of this work is twofold. Firstly, as already mentioned, an arbitrary finite reflection group is admitted, hence full generality 
(with respect to this kind of symmetry) is reached. It is worth pointing out that the case of a finite number of orthogonal mirrors considered in \cite{MS} 
carried the character of 'product case'. For an arbitrary finite reflection group in general there is no orthogonality between the mirrors and, in some sense, 
this is the heart of the matter. Clearly, treatment of an arbitrary reflection group required completely new arguments. This is exhibited in using weighted extension and averaging operators as well as some additional technical tools. See Sections \ref{sec:prel} and \ref{sec:aux} for details. Generality mentioned 
above permits, for instance, to obtain closed formulas for the Neumann and Dirichlet heat kernels on cones on the plane with apertures $\pi/n$, $n\ge3$; 
up to our knowledge these formulas are new for $n\neq 2^j$. The same can be said about analogous truncated cones (intersections of these infinite cones with the 
unit disc centered at the origin) on the plane, and also about polyhedral cones in $\mathbb R^d$ being positive Weyl chambers corresponding to a finite reflection  group 
and their truncated versions.

Secondly, the new concept of using homomorphisms from ${\rm Hom}(W,\widehat{\mathbb Z}_2)$ in construction of self-adjoint extensions of 
$-\Delta_{\Omp}$ is presented. Heat kernels associated with these extensions are useful tools in solving mixed Neumann-Dirichlet initial-boundary value
problems on cones being positive Weyl chambers. This issue is outlined in Section \ref{sec:ortho} and deserves further studies. 

Finally we mention one more aspect of novelty. Namely, Theorem \ref{thm:main} allows one to obtain heat kernel formulas for sets with less regular boundary 
from heat kernel formulas for more regular sets. Often a general theory applies to sets $\Om$ with, say, $C^{1,1}$ boundary while the set $\Omp$ obtained 
by cutting off from such $\Om$ the positive part, $\Omp=\Om\cap C_+$,  may have only Lipschitz boundary.

The paper is organized as follows. Section \ref{sec:prel} contains definitions of weigh\-ted extension and averaging operators and facts about them needed subsequently. Section \ref{sec:aux} is devoted to the proofs of auxiliary results which are used in the proof of the main result in Section \ref{sec:proof}. 
Section \ref{sec:ext} contains definition and basic results on $\eta$-Laplacians, which are self-adjoint extensions in $L^2(\Omp)$ of the Laplacian $-\Delta_{\Omp}$.
Finally, in Section \ref{sec:app} we discuss two basic examples, the case of a finite reflection group associated with an orthonormal root system and the case of dihedral group.  

\section{Preliminaries} \label{sec:prel}
From now on till Section \ref{sec:proof}, $\Om$ is a fixed open subset of $\R$, symmetric with respect to a given finite reflection group $W$ associated with a fixed root system $R$. Consequently, $\Omp$ will denote its positive part related to a chosen (and fixed) subsystem $R_+$ of positive roots. 
Since the case $d=1$ is simple (and is contained in \cite{MS}), henceforth we assume that $d\ge 2$. Throughout this and the next sections we obey the following convention: we use small letters (like $\phi$, $f$,\ldots) to denote functions on $\Omp$, and capital letters (like $\Phi$, $F$,\ldots) to denote functions on 
$\Om$. The symbol $\v$ is reserved to denote a compactly supported $C^\infty$ function with support either in $\Om$ or $\Omp$ (this always will be clear from the context). For $g\in W$ and a function $F$ on $\Om$ by $F_g$ we mean the function $F_g(x):=F(gx)$, $x\in\Om$. Also, recall that the norm 
$\|\cdot\|_{H^1(\Om)}$ in $H^1(\Om)$ is given by
$$
\|F\|_{H^1(\Om)}=\Big( \|F\|_{L^2(\Om)}^2+\sum_{i=1}^d \|\partial_i F\|_{L^2(\Om)}^2\Big)^{1/2},
$$
where $\partial_i F$ are weak partial derivatives of $F$ on $\Om$. Similarly for $\|\cdot\|_{H^1(\Omp)}$ in  $H^1(\Omp)$.

Let $\omega\colon W\to \C$ be a function on $W$, called a\textit{ weight} henceforth. Recall that $W$ acts on the set of 
Weyl chambers simply transitively.
\begin{definition}
The \textit{weighted extension operator} $\EE^\omega$ acts on  functions $\phi$ defined on $\Omp$ by
$$
\EE^\omega \phi(gx)=\omega(g)\phi(x), \qquad x\in\Omp,\quad g\in W.
$$
The \textit{weighted averaging operator} $\A_\omega$ acts on  functions $\Phi$ defined on $\Om$ by
$$
\A_\omega \Phi(y)=\frac1{|W|}\sum_{g\in W}\omega(g)\Phi(gy), \qquad y\in\Om.
$$
\end{definition}
Clearly, the resulting functions $\EE^\omega \phi$ and $\A_\omega \Phi$ live on $\Om$ (to be precise, $\EE^\omega \phi$ is defined on $\Om$ up to a set of 
Lebesgue measure zero). Also, for a function $\Phi$ on $\Om$, by $\Phi^+$ we will mean the restriction of $\Phi$ to $\Omp$, $\Phi^+=\Phi|_{\Omp}$.

For the weight $\omega\equiv 1$ we shall simply write $\EE$ and $\A$ rather than $\EE^1$ and $\A_1$. For the weight ${\rm sgn}$ we write $\EE^{{\rm sgn}}$ 
and $\A_{{\rm sgn}}$. Then $\EE\phi$ and $\EE^{{\rm sgn}}\phi$ can be called the \textit{even} and \textit{odd extensions} of $\phi$ onto $\Om$, respectively.

It is easily seen that the following averaging invariance properties hold
$$
\A(\EE\phi)=\EE\phi, \qquad \A_{{\rm sgn}}(\EE^{{\rm sgn}}\phi)=\EE^{{\rm sgn}}\phi.
$$
\begin{lemma} \label{lem:first}
For suitable functions $\phi$ and $\Phi$ on $\Omp$ and $\Om$ (i.e. such that the first integral below makes sense), respectively, we have
\begin{equation}\label{H1}
\int_\Om\EE^\omega \phi\cdot \Phi=|W|\int_{\Om_+}\phi\cdot (\A_\omega \Phi)^+.
\end{equation}
\end{lemma}
\begin{proof}
We write, using the change of variables $y:=gx$, $y\in g\Omp$, $x\in\Omp$,
\begin{align*}
\int_\Om\EE^\omega \phi\cdot \Phi
&=\sum_{g\in W}\int_{g\Omp}\EE^\omega \phi(y)\cdot \Phi(y)\,dy\\
&=\sum_{g\in W}\int_{\Omp} \omega(g)\phi(x)\cdot \Phi(gx)\,dx\\
&=\int_{\Omp} \phi(x)\cdot \big(\sum_{g\in W}\omega(g)\Phi(gx)\big)\,dx\\
&=|W|\int_{\Omp} \phi\cdot \big(\A_\omega\Phi\big)^+.
\end{align*}
\end{proof}

For $g\in W$ let $[g_{ij}]$ be the matrix of the linear transformation $x\mapsto gx$. This means that for $x=(x_1,\ldots, x_d)$ and 
$gx=((gx)_1,\ldots,(gx)_d)$,
$$
(gx)_i=\sum_j g_{ij}x_j, \qquad i=1,\ldots,d.
$$
Clearly, $[g_{ij}]$ is also the Jacobi matrix of this transformation,
$$
g_{ij}=\partial_j((gx)_i).
$$
Moreover, since $[g_{ij}]$ is an orthogonal matrix, we have the orthogonality of rows and columns, $\sum_j g_{ij}g_{mj}=\delta_{im}$ and 
$\sum_i g_{ij}g_{im}=\delta_{jm}$.

\begin{lemma} \label{lem:firstbis}
The operator $\A_\omega$ acts boundedly on $H^1(\Om)$. Moreover,  $\A_\omega$ leaves $H_0^1(\Om)$ invariant.
\end{lemma}
\begin{proof}
If  $F\in H^1(\Om)$, then it is immediately seen that for any $g\in W$ the function $F_g$ is also in $H^1(\Om)$ and, since 
$(\partial_i F_g)(x)=\sum_{j=1}^d\partial_jF(gx)g_{ji}$, the mapping $F\to F_g$ is bounded. (See \cite[Theorem 3.41]{AF} for a general result on involving coordinate transformation in Sobolev space context.) Consequently, $\A_\omega$ maps $H^1(\Om)$ into itself and is bounded. The last claim in the statement 
of the lemma is obvious.
\end{proof}

In what follows, for any given $j,i\in\{1,\ldots,d\}$, by writing $\A_{g_{ji}}$ or $\A_{\om\cdot g_{ji}}$  we shall mean the weighted averaging  operators 
with the weight functions $g\to g_{ji}$ or $g\to \om(g) g_{ji}$ on $W$, respectively. Similarly for $\EE^{g_{ji}}$ and $\EE^{\om\cdot g_{ji}}$.
\begin{lemma} \label{lem:second}
For any  $F\in H^1(\Om)$ we have
\begin{equation}\label{H2}
\A_\om(\partial_j F)=\sum_{i=1}^d \partial_i(\A_{\om\cdot g_{ji}}F)
\end{equation}
and 
\begin{equation}\label{H3}
\partial_i (\A_\om F)=\sum_{j=1}^d \A_{\om\cdot g_{ji}}(\partial_j F).
\end{equation}
\end{lemma}
\begin{proof}
First, assume $F\in C^1(\Om)$. By the very definition
$$
\A_{\om\cdot g_{ji}} F(y)=\frac1{|W|}\sum_{g\in W} \om(g) g_{ji} F(gy), \qquad y\in\Om,
$$
and hence 
\begin{align*}
\partial_i(\A_{\om\cdot g_{ji}} F)(y)&=\frac1{|W|}\sum_{g\in W} \om(g)g_{ji}\partial_i\big( F(gy)\big)\\
&=\frac1{|W|}\sum_{g} \om(g)g_{ji}\sum_{m}\partial_m F(gy)\cdot g_{mi}.
\end{align*}
For a given $j$, summing over $i$ gives \eqref{H2}. 

For a general $F\in H^1(\Om)$ find a sequence $\{F_n\}\subset C^1(\Om)\cap H^1(\Om)$ such that $F_n\to F$ in $H^1(\Om)$ (see \cite[Theorem 3.17]{AF}). 
This means, in particular, that $F_n\to F$ in  $L^2(\Om)$, and since $\A_\om$ is bounded on $L^2(\Om)$, also $\A_\om F_n\to \A_\om F$ in $L^2(\Om)$. 
On the other hand, $\A_{\om\cdot g_{ji}}$ is bounded on $H^1(\Om)$, so $\A_{\om\cdot g_{ji}} F_n\to \A_{\om\cdot g_{ji}}F$ in $H^1(\Om)$, and hence 
$\partial_i(\A_{\om\cdot g_{ji}} F_n)\to \partial_i(\A_{\om\cdot g_{ji}}F)$ in $L^2(\Om)$, as $n\to \infty$. Since \eqref{H2} holds for any such $F_n$, consequently, \eqref{H2} follows.

The verification of \eqref{H3} is analogous.
\end{proof}

\section{Auxiliary results} \label{sec:aux}

In the sequel we shall need the explicit description of the boundary of $\overline{C_+}$, the closure of the positive Weyl chamber $C_+$ distinguished 
by a set of positive roots $R_+$. For this purpose we recall the notion of the \textsl{system of simple roots}; see \cite{Hu}, \cite{DX} or \cite{He}. 
This is the unique subset $\Pi=\{\a_1,\ldots,\a_m\}$ of $R_+$ which is a basis of ${\rm lin}\{\a\colon \a\in R_+\}$  (so $1\le m\le d$)  and each $\a\in R_+$ 
is a linear combination of $\a_1,\ldots,\a_m$ with non-negative coefficients. Consequently, 
$$
C_+=\{x\in\Om\colon \langle x,\a_k\rangle>0,\,\, k=1,\ldots,m\},
$$ 
and the closure $\overline{C_+}$ has exactly $m$ facets $\overline{C_+}\cap \langle\a_k\rangle^\bot$, $k=1,\ldots,m$. These facets are closed $(d-1)$-dimensional infinite cones in the hyperplanes $\langle\a_k\rangle^\bot$, respectively. If $m=1$, then this single cone coincides with $\langle\a_1\rangle^\bot$, 
otherwise, for $m\ge 2$, these cones are proper. Finally we mention that the 'simple' reflections $\sigma_{\a_k}$, $k=1,\ldots,m$, generate $W$. 

In what follows, if $V$ is a linear space of functions on $\Om$, then we shall write $V^+$ for the linear space of the restrictions of functions from $V$ 
to $\Omp$, $V^+:=\{\Phi^+\colon \Phi\in V\}$. Analogously, by $\A_\omega V$ we shall mean the space of (weighted) averages of functions from $V$, 
$\A_\omega V:=\{\A_\omega\Phi\colon \Phi\in V\}$. 

\begin{proposition} \label{pro:ext}
Let $f\in H^1(\Om_+)$. Then $\EE f\in H^1(\Om)$ and, 
\begin{equation}\label{bb}
\partial_j(\EE f)=\sum_{i=1}^d\EE^{g_{ji}}(\partial_i f).
\end{equation}
\end{proposition}
\begin{proof}
If $f\in C^1_c(\Omp)$, then obviously $\EE f\in C^1_c(\Om)$ (we set $\EE f(x)=0$ for $x\in\Om\cap \bigcup_{\a\in R}\langle \a\rangle^\bot$ 
provided this set is nonempty) and checking that 
$$
\partial_j(\EE f)(gx)=\sum_{i=1}^dg_{ji}\partial_i f(x),\qquad x\in\Omp,\quad g\in W,
$$
which is just the explicit form of \eqref{bb}, relies on a direct calculation.

In the case of general $f\in H^1(\Om_+)$ we shall verify that for any  $ \v\in C_c^\infty(\Om)$ it holds
\begin{equation}\label{aa}
\int_{\Om_+}f\, \big(\sum_i\partial_i(\A_{g_{ji}}  \v)\big)^+=-\int_{\Om_+}\sum_i\partial_i f\, (\A_{g_{ji}}  \v)^+.
\end{equation}
This will be sufficient for drawing the conclusion  that $\partial_j(\EE f)$ exist and \eqref{bb} holds true, which also means that $\EE f\in H^1(\Om)$. 
Indeed, to check \eqref{bb} with the aid of \eqref{aa} we use successively Lemma \ref{lem:first}, \eqref{H2} for $\om\equiv1$, \eqref{aa} and  
Lemma \ref{lem:first} again, to write
\begin{align*}
\int_\Om \EE f\, \partial_j\v=|W|\int_{\Omp} f\, \big(\A(\partial_j\v)\big)^+&=|W|\int_{\Om_+} f\, \big(\sum_i\partial_i(\A_{g_{ji}}\v)\big)^+\\
&=-|W|\int_{\Om_+}\sum_i\partial_i f\, (\A_{g_{ji}}\v)^+\\
&=-\int_\Om \big(\sum_i\EE^{g_{ji}}(\partial_i f)\big)\,\v.
\end{align*}

The remaining part of the proof is devoted to verification of  \eqref{aa}. Assume $f\in H^1(\Om_+)$ and $\v\in C_c^\infty(\Om)$ are fixed and recall that  
$\overline{C_+}$ has exactly $m$ facets: $\overline{C_+}\cap \langle\a_k\rangle^\bot$, $k=1,\ldots,m$. We now choose $\widetilde{\Om}$, a bounded open 
set with $C^1$ boundary (piecewise $C^1$ would be sufficient as well), which is $W$-\,symmetric and such that ${\rm supp}\,\v\subset \widetilde{\Om}\subset \Om$. 
Observe that the assumption on support of $\v$ and symmetry of $\widetilde{\Om}$ imply that also ${\rm supp}\,\A_\omega\v\subset \widetilde{\Om}$ for any 
weight $\omega$. Observe also, that verification of \eqref{aa} reduces to checking that for any fixed $j\in \{1,\ldots,d\}$ and 
$\widetilde{\Om}_+:=\widetilde{\Om}\cap C_+$ (recall that $\Omp=\Om\cap C_+$) it holds
\begin{equation}\label{cc}
\int_{\widetilde{\Om}_+}f\, \big(\sum_i\partial_i(\A_{g_{ji}}  \v)\big)^+=-\int_{\widetilde{\Om}_+}\sum_i\partial_i f\, (\A_{g_{ji}}  \v)^+.
\end{equation}

Since $\widetilde{\Om}_+$ is a bounded  open set with piecewise $C^1$ boundary, a general theory applies and for any $i\in \{1,\ldots,d\}$, the following integration by parts formula holds, 
\begin{equation}\label{dd}
\int_{\widetilde{\Om}_+}g\,\partial_i h=-\int_{\widetilde{\Om}_+}\partial_i g\, h+\int_{\partial(\widetilde{\Om}_+)}gh\,\nu_i\, d\s, 
\qquad g,h\in H^1(\widetilde{\Om}_+).
\end{equation}
Here $\nu=\nu(x)=(\nu_1,\ldots,\nu_d)$ denotes the outward unit normal vector at $x\in \partial(\widetilde{\Om}_+)$ (whenever it exists) and $d\s$ means 
the surface measure on $\partial(\widetilde{\Om}_+)$, the boundary of $\widetilde{\Om}_+$. See \cite[pp. 263-271]{A}, where \eqref{dd} is stated with 
weaker assumption for bounded sets with Lipschitz boundary, or \cite[Theorem D.8]{Sch}, where assumptions on  boundary are stronger. 
It should be also pointed out that in \eqref{dd} the boundary values of $g$ and $h$ (\textit{traces} of  $g$ and $h$) are well defined 
$L^2(\partial(\widetilde{\Om}_+),\, d\s)$ functions; this is the content of the \textit{trace theorem}, see \cite[ A8.6]{A} or \cite[Theorem D.6]{Sch}. 

Now, substituting in \eqref{dd} $f$ for $g$ and $\A_{g_{ji}}  \v$ for $h$, and summing over $i$, we obtain an identity like \eqref{cc} but with the additional term
\begin{equation}\label{zz}
\int_{\partial(\widetilde{\Om}_+)}f\,\big(\sum_i\A_{g_{ji}}\v\cdot\nu_i\big)\, d\s
\end{equation}
on the right hand side. Aiming to be pedantic we should keep in mind that in fact we substitute  in \eqref{dd} the restrictions $f|_{\widetilde{\Om}_+}$ 
and $(\A_{g_{ji}}\v)|_{\widetilde{\Om}_+}$ for $g$ and $h$, respectively. But both restrictions are indeed in $ H^1(\widetilde{\Om}_+)$ and, moreover, 
$\partial_i\big(f|_{\widetilde{\Om}_+}\big)=\big(\partial_i f \big)|_{\widetilde{\Om}_+}$ and 
$\partial_i\big((\A_{g_{ji}}\v)|_{\widetilde{\Om}_+}\big)=\big(\partial_i( \A_{g_{ji}}\v) \big)|_{\widetilde{\Om}_+}$. Also, $f$ in \eqref{zz} is understood 
as the trace of $f|_{\widetilde{\Om}_+}$ onto $\partial(\widetilde{\Om}_+)$. The trace exists since $f|_{\widetilde{\Om}_+}$ is in $H^1(\widetilde{\Om}_+)$ 
and $\widetilde{\Om}_+$ is bounded with Lipschitz boundary.

To conclude verification of  \eqref{cc} and thus to finish the proof of the proposition we will check that the quantity in \eqref{zz} vanishes.

The boundary $\partial(\widetilde{\Om}_+)$ consists of at most $m$ 'flat' parts, $\overline{\widetilde{\Om}_+}\cap \langle\a_k\rangle^\bot$, $ k\in\{1,\ldots,m\}$,  
(some of these sets or even all of them may be empty) and the 'irregular' part $\partial \widetilde{\Om}\cap C_+$. 
On the latter part all $\A_{g_{ji}}\v$, $i=1,\ldots,d$, vanish and thus the corresponding part of integration over the boundary vanishes. In fact, also on each nonempty 'flat' part of the boundary the expression $\sum_i\A_{g_{ji}}\v\cdot\nu_i$ vanishes as well. Indeed,  note that $\nu$,  the outward unit normal vector to the 'flat' part $\overline{\widetilde{\Om}_+}\cap \langle\a_k\rangle^\bot$, 
which we assume to be nonempty, coincides with $\a_k$. Let 
$\a_k=(\a_{k1},\ldots, \a_{kd})$. Our aim is to check (slightly more than required, with $\v$ extended outside $\Om$ by zero) that
\begin{equation}\label{ee}
\sum_i\A_{g_{ji}}\v(x)\,\a_{ki}=0, \qquad x\in \langle\a_k\rangle^\bot.
\end{equation}
Let $W^+={\rm ker(sgn)}=\{g\in W\colon {\rm sgn}(g)=1\}$. Then, due to the partition $W=W^+\sqcup (W^+\cdot\s_{\a_k})$, we have for $x\in \langle\a_k\rangle^\bot$
\begin{align*}
|W|\sum_i \A_{g_{ji}}\v(x)\,\a_{ki}&=\sum_i \big(\sum_{g\in W} g_{ji}\v(gx) \big)\, \a_{ki}\\
&=\sum_i \Big(\sum_{g\in W^+} g_{ji}\v(gx)+ \sum_{g\in W^+} (g\s_{\a_k})_{ji}\v(g\s_{\a_k}x)\Big)\, \a_{ki}\\
&=\sum_i \Big(\sum_{g\in W^+}\v(gx) \big(g_{ji}+ (g\s_{\a_k})_{ji}\big)\Big)\, \a_{ki}\\
&=\sum_{g\in W^+}\v(gx)\sum_i \big(g_{ji}+ (g\s_{\a_k})_{ji}\big)\, \a_{ki} .
\end{align*}
Since $\v$ is arbitrary, to reach \eqref{ee}, our goal, we should expect the innermost sum in the  above last line to vanish. This is indeed the case. For $g\in W^+$ 
(and  $ k\in\{1,\ldots,m\}$ fixed) an easy matrix calculation shows that 
$$
(g\s_{\a_k})_{ji}=g_{ji}-2\a_{ki}\sum_sg_{js}\a_{ks},
$$
and hence 
\begin{equation*}
\sum_i \big(g_{ji}+ (g\s_{\a_k})_{ji}\big)\, \a_{ki}=2\Big(\sum_i g_{ji}\a_{ki}  - \sum_i \a_{ki}^2\cdot \sum_sg_{js}\a_{ks}\Big).
\end{equation*}
The last expression vanishes since $\a_k=(\a_{k1},\ldots,\a_{kd})$ is a unit vector.

This finishes checking of \eqref{ee}, concludes verification of \eqref{cc} and thus also \eqref{aa}, and hence finishes the proof of Proposition \ref{pro:ext}.
\end{proof}

The following result is an immediate consequence of the previous one.

\begin{proposition} \label{pro:Sob1}
We have
\begin{equation}\label{S1}
H^1(\Om_+)=\A(H^1(\Om))^+. 
\end{equation}
\end{proposition}
\begin{proof}
The inclusion $\supset$ in \eqref{S1} is obvious: if $F\in H^1(\Om)$, then, by Lemma \ref{lem:firstbis}, $\A F\in H^1(\Om)$ and hence $(\A F)^+\in H^1(\Om_+)$.
To prove the opposite inclusion $\subset$  take $f\in H^1(\Om_+)$.  By Proposition \ref{pro:ext} we have $\EE f\in H^1(\Om)$, and hence  $f=(\EE f)^+=\big(\A(\EE f)\big)^+$ is in $\A(H^1(\Om))^+$.
\end{proof}
Clearly, for an arbitrary weight the result of Proposition \ref{pro:ext} is not true. However, with a stronger assumption we have the following.

\begin{proposition} \label{pro:extbis}
Let $f\in H^1_0(\Om_+)$ and $\om$ be a weight on $W$. Then $\EE^\om  f\in H^1_0(\Om)$ and
\begin{equation}\label{bbbis}
\partial_j(\EE^\om  f)=\sum_{i=1}^d\EE^{\om\cdot g_{ji}}(\partial_i f).
\end{equation}
\end{proposition}
\begin{proof}
If $f\in C^1_c(\Omp)$, then obviously $\EE^\om  f\in C^1_c(\Om)$  and checking that 
$$
\partial_j(\EE^\om  f)(gx)=\sum_{i=1}^d\om (g) g_{ji}\,\partial_i f(x),\qquad x\in\Omp,\quad g\in W,
$$
which is just the explicit form of \eqref{bbbis}, relies on a direct calculation.

In the case of general $f\in H^1_0(\Om_+)$ we shall verify that for any weight $\om'$ and $ \v\in C_c^\infty(\Om)$ it holds
\begin{equation}\label{aabis}
\int_{\Om_+}f\, \partial_i(\A_{\om'}  \v)^+=-\int_{\Om_+}\partial_i f\, (\A_{\om'}  \v)^+, \qquad  i=1,\ldots,d.
\end{equation}
This will be sufficient for drawing the conclusion that $\partial_j(\EE^\om   f)$ exist and, moreover, \eqref{bbbis} holds true, which also means that 
$\EE^\om f\in H^1(\Om)$. Indeed, to check \eqref{bbbis} with the aid of \eqref{aabis}, we use Lemma \ref{lem:first} twice, \eqref{H2}, and \eqref{aabis} 
applied to $\om'= \om \cdot g_{ji}$, to write
\begin{align*}
\int_\Om \EE^\om  f\,\partial_j\v=|W|\int_{\Omp} f\, \big(\A_\om (\partial_j\v)\big)^+&=|W|\int_{\Om_+} f\, \big(\sum_i\partial_i(\A_{\om \cdot g_{ji}}\v)\big)^+\\
&=-|W|\sum_i\int_{\Om_+}\partial_i f\, (\A_{\om \cdot g_{ji}}\v)^+\\
&=-\int_\Om \big(\sum_i\EE^{\om \cdot g_{ji}}(\partial_i f)\big)\,\v.
\end{align*}

We now focus on the proof  of  \eqref{aabis}. Notice that it suffices to verify this identity for $f\in C^\infty_c(\Omp)$. Indeed, assuming \eqref{aabis} 
holds true in such setting, take $\v_n\to f$ in $H^1(\Om_+)$, $\v_n\in C^\infty_c(\Omp)$, which means that, in particular, $\v_n\to f$ and $\partial_i\v_n\to \partial_if$ in $L^2(\Om_+)$ as $n\to\infty$. Thus, considering  \eqref{aabis} with these $\v_n$'s in place of $f$, and letting $n\to\infty$ gives  
\eqref{aabis} in its original form.

Assume therefore that $f\in C^\infty_c(\Om_+)$ and $\v\in C_c^\infty(\Om)$ are fixed, choose $\Om_0\subset \Om_+$, a bounded open set with $C^1$ boundary, 
such that ${\rm supp}\,f\subset \Om_0$, and observe that verification of \eqref{aabis} reduces to checking that  
\begin{equation}\label{ccbis}
\int_{\Om_0}f\, \partial_i(\A_{\om'}  \v)^+=-\int_{\Om_0}\partial_i f\, (\A_{\om'}  \v)^+.
\end{equation}
With this choice of  $\Om_0$ and by using the integration by parts formula (for smooth functions), checking \eqref{ccbis} further reduces to verification 
that the quantity as in \eqref{zz} but with $\A_{g_{ji}}$ replaced by $\A_{\om'}$ and $\widetilde{\Om}_+$ replaced by $\Om_0$, vanishes. But this is 
obvious since $f$ vanishes on $\partial\,\Om_0$.

Finally, we show that in fact, for $f\in H^1_0(\Om_+)$ we have  $\EE^\om   f\in H^1_0(\Om)$. Indeed take $\{\v_n\}\subset C^\infty_c(\Omp)$ such that $\v_n\to f$ in $H^1_0(\Om_+)$. Since $\{\v_n\}$ is a 
Cauchy sequence in $H^1(\Om_+)$, it follows that $\{\EE^{\om }\v_n\}$ is a Cauchy sequence in $H^1(\Om)$. (Clearly $\EE^{\om }(\v_n-\v_m)\to0$ in $L^2(\Om)$ as $n,m\to\infty$, so verification 
that for every $j\in\{1,\ldots,d\}$, also $\partial_j(\EE^{\om }(\v_n-\v_m))\to\infty$ in $L^2(\Om)$ reduces, by using \eqref{bbbis}, to noting that for every $i\in\{1,\ldots,d\}$, 
$\EE^{\om \cdot g_{ji}}(\partial_j(\v_n-\v_m))\to0$ in $L^2(\Om)$.) 
Let $F=\lim_{n\to\infty}\EE^{\om }\v_n$ in $H^1(\Om)$. Obviously, $\{\EE^{\om }\v_n\}\subset C^\infty_c(\Om)$ and hence $F\in H^1_0(\Om)$ and $F|_{\Omp}=f$.
\end{proof}

\begin{proposition} \label{pro:Sob2}
We have
\begin{equation}\label{S2}
 H^1_0(\Om_+)=\A_{\rm sgn}(H^1_0(\Om))^+.
\end{equation}
\end{proposition}
\begin{proof}
We begin with the inclusion $\subset$ in \eqref{S2} and take $f\in H^1_0(\Om_+)$. Proposition \ref{pro:extbis} says that $\EE^{\rm sgn} f\in H^1_0(\Om)$, 
and hence  $f=(\EE^{\rm sgn} f)^+=\big(\A_{\rm sgn}(\EE^{\rm sgn} f)\big)^+$ is in $\A_{\rm sgn}(H^1_0(\Om))^+$.

Proving the opposite inclusion $\supset$ take $F\in H^1_0(\Om)$. Then, by Lemma \ref{lem:firstbis},  $\A_{\rm sgn}F\in H^1_0(\Om)$ and hence 
$(\A_{\rm sgn}F)^+\in H^1(\Omp)$. We will show that, in fact, $(\A_{\rm sgn}F)^+\in H^1_0(\Omp)$. Let $\{\Phi_n\}\subset C_c^\infty(\Om)$ be such 
that $\Phi_n\to F$ in $ H^1(\Om)$. Then, again by Lemma \ref{lem:firstbis},  also $\A_{\rm sgn}\Phi_n\to\A_{\rm sgn}F$ in $H^1(\Om)$ 
and hence $(\A_{\rm sgn}\Phi_n)^+\to(\A_{\rm sgn}F)^+$ in $ H^1(\Omp)$. Since $\{\A_{\rm sgn}\Phi_n\}\subset C_c^\infty(\Om)$ and each $\A_{\rm sgn}\Phi_n$  
for every $\a\in R$ satisfies $\A_{\rm sgn}\Phi_n(\s_\a y)=-\A_{\rm sgn}\Phi_n(y)$, $y\in\Om$, it suffices to verify that each restriction $(\A_{\rm sgn}\Phi_n)^+$ can be approximated in the $H^1(\Omp)$ norm by $C^\infty_c(\Omp)$ functions. 

In fact we shall prove slightly more general claim. Namely, assuming $\Phi\in C^\infty_c(\Om)$ is such that  $\Phi$ vanishes on $\langle \a_k\rangle^\bot$ for every $k=1,\ldots,m$,
we shall verify that $\Phi^+$ can be approximated in $H^1(\Omp)$ by $C^\infty_c(\Omp)$ functions. 

Let $\b\in C^\infty(0,\infty)$ be such that $\b(t)=0$ for $0<t<1/2$ and $\b(t)=1$ for $t>1$ 
and $\|\b\|_\infty=1$. We define the sequence of functions $\widetilde{\b_N}$ on $\Omp$, $N\in\N$, by
$$
\widetilde{\b_N}(x)=\Pi_{k=1}^m\b(N\delta_k(x)), \qquad x\in \Omp,
$$
where $\delta_k(x)=\sum_{i=1}^d\a_{ki}x_i$, $\a_k=(\a_{k1},\ldots,\a_{kd})$, denotes the $\ell^1$ distance of $x\in\Omp$ to the hyperplane 
$\langle \a_k\rangle^\bot$. Consequently, we define 
$$
\widetilde{\Phi_N}(x)=\widetilde{\b_N}(x)\Phi(x), \qquad x\in \Omp.
$$
Clearly, $\widetilde{\Phi_N}\in C^\infty_c(\Omp)$. We shall verify that $\widetilde{\Phi_N}\to \Phi^+$ in $H^1(\Omp)$, as $N\to\infty$. We have
$$
\|\widetilde{\Phi_N}-\Phi^+\|_{H^1(\Omp)}^2=\|\Phi^+(\widetilde{\b_N}-1)\|_{L^2(\Omp)}^2+\sum_{j=1}^d \|\partial_j(\Phi^+(\widetilde{\b_N}-1))\|_{L^2(\Omp)}^2.
$$
Obviously, the first summand on the right hand side goes to 0 as $N\to \infty$. It remains to check that the same happens for each of the remaining $d$ summands. 
We have
$$
\partial_j(\Phi^+(\widetilde{\b_N}-1))=\partial_j\Phi^+(\widetilde{\b_N}-1)+\Phi^+\partial_j\widetilde{\b_N},
$$
so we need to consider the second summand only (note that $\partial_j\Phi^+=(\partial_j\Phi)^+$). For simplicity, consider $j=1$. Then
$$
\partial_1\widetilde{\b_N}(x)=N\sum_{k=1}^m \partial_1\delta_k(x)\b'(N\delta_k(x))\prod_{s\neq k}\b(N\delta_s(x)).
$$
Note that $\partial_1\delta_k(x)$ is a constant not depending on $x\in\Omp$. We are therefore reduced to show that for any $k\in\{1,\ldots,m\}$ the quantity
\begin{equation}\label{max}
N^2\int_{\Omp}|\Phi^+(x)\b'(N\delta_k(x))\prod_{s\neq k}\b(N\delta_s(x))|^2\,dx
\end{equation}
tends to 0 as $N\to\infty$. Assume that ${\rm supp}\,\Phi\,\subset B_R(0)$ with some $R>0$ and let $\|\b'\|_\infty=C$. There exists a constant $M$ such that for $x\in\Omp$ we have
$$
|\Phi(x)|\le M\min_{1\le l\le m}\delta_l(x).
$$
This follows by the mean value theorem and the fact that $\Phi$ is compactly supported and vanishes on each facet $\langle \a_l\rangle^\bot$. Therefore \eqref{max} is estimated from above by
$$
N^2C^2M^2\int_{A(N,R)}\big(\min_{1\le l\le m}\delta_l(x)\big)^2dx,
$$
where $A(N,R)=B_R(0)\cap \{x\in C_+\colon \min_{1\le l\le m}\delta_l(x)\le 1/N\}$. It is now clear that the last quantity tends to 0 as $N\to\infty$.

This finishes the proof of the claim and thus concludes the proof of the $\supset$. 
\end{proof}
At the end of this section it is worth noting that Propositions \ref{pro:ext} and \ref{pro:extbis} can be regarded as \textit{extension theorems} in the sense of \cite[Section 5]{AF}.
Although it is not contained in the statements of the propositions, in fact it follows from the proofs that both extension operators are bounded, that is $\|\EE f\|_{H^1(\Om)}\lesssim\|f\|_{H^1(\Omp)}$, 
and similarly for $\EE^{\omega}$.

\section{$\eta$-extensions of the Laplacian} \label{sec:ext}
Let $\eta\colon W\to \widehat{\mathbb Z}_2$ be a homomorphism. For a function $F$ on $\Om$ we say it is $\eta$\,- symmetric, provided for every $g\in W$ and $y\in\Om$ it holds $F(gy)=\eta(g)F(y)$. We call the averaging operator $\A_\eta$ the $\eta$-symmetrization operator. Obviously, $F$ is 
$\eta$\,- symmetric if and only if $\A_\eta F=F$. Clearly, $\A_\eta F$ and $\EE^\eta \phi$ are $\eta$-symmetric whatever $F$ and $\phi$ are, and hence 
$$
\A_\eta (\A_\eta F)=\A_\eta F\quad {\rm and}\quad \A_\eta (\EE^\eta \phi)=\EE^\eta \phi.
$$ 
Recall that for any given $i,j\in\{1,\ldots,d\}$, by $\A_{\eta\cdot g_{ji}}$  we mean the weighted averaging  operator with the 
weight function $g\to \eta(g) g_{ij}$ on $W$. Similarly for  $\EE^{\eta\cdot g_{ij}}$.
Also, if $V$ is a linear space of functions on $\Om$, then we write $V^+$ for the linear space of the restrictions of functions from $V$ to $\Om_+$,
$V^+:=\{\Phi^+\colon \Phi\in V\}$, and $\A_\eta V$ for the linear space of $\eta$-symmetrizations of functions from $V$, $\A_\eta V:=\{\A_\eta\Phi\colon \Phi\in V\}$.

It will be convenient to use the following notation
$$
H^1_\eta(\Om_+):=\big(\A_\eta H^1(\Om)\big)^+.
$$
Obviously, since $\A_\eta$ maps $H^1(\Om)$ into itself, $H^1_\eta(\Om_+)$ is a subspace in $H^1(\Om_+)$. It is easily seen that $H^1_\eta(\Om_+)$ should be regarded as the space 
of restrictions to $\Om_+$ of $\eta$-symmetric functions from $H^1(\Om)$, or as the space of functions from $L^2(\Om_+)$, with $\eta$-symmetric extensions in $H^1(\Om)$,
$$
H^1_\eta(\Om_+)=\{F\in H^1(\Om)\colon F=\A_\eta F\}^+=\{f\in L^2(\Om_+)\colon \EEeta f\in H^1(\Om)\}.
$$

From now on till the end of this section, in the statements of results, $\eta\in{\rm Hom}(W,\,\widehat{\mathbb Z}_2)$ is a fixed homomorphism. 
The following result is similar to Propositions \ref{pro:ext} and \ref{pro:extbis}. Note, however, that the present assumptions are much stronger and hence the reasoning in the proof is much simpler. 
\begin{proposition} \label{lem:ext}
Let $f\in H^1_\eta(\Om_+)$. Then 
\begin{equation}\label{bbb}
\partial_j(\EEeta f)=\sum_{i=1}^d\EE^{\eta\cdot g_{ji}}(\partial_i f).
\end{equation}
\end{proposition}
\begin{proof}
Take $f\in H^1_\eta(\Om_+)$. This means that $\EEeta f\in H^1(\Om)$ and, in particular, $f\in H^1(\Om_+)$, and   hence  both sides of \eqref{bbb} make sense. 
To verify  \eqref{bbb} it suffices to check that both sides coincide on each Weyl chamber $g\Om_+$, $g\in W$. But
$$
\partial_j(\EEeta f)|_{g\Om_+}=\partial_j\big((\EEeta f)|_{g\Om_+}\big)
$$
and hence the task reduces to checking that the $j$th weak partial derivative of the function $gx\mapsto \eta(g)f(x)$ defined on  $g\Om_+$, $x\in \Om_+$,
coincides with the restriction to $g\Om_+$ of the right-hand side of \eqref{bbb}, i.e. the function  $gx\mapsto \eta(g)\sum_{i=1}^dg_{ji}\partial_i f(x)$.
This follows by an easy calculation.
\end{proof}

The following result will be crucial in the definition of $-\Delta^\eta_+$, the $\eta$-Laplacian  on $\Om_+$.
\begin{proposition} \label{pro:clo}
$H^1_\eta( \Om_+)$ is a closed subspace in $H^1(\Om_+)$.
\end{proposition}
\begin{proof}
 We begin with  remark that  if  $h\in L^2(\Omp)$  and $a\in\mathbb R$, then 
\begin{equation}\label{eq}
\|\EE^{\eta\cdot a} h\|_{L^2(\Om)}=|W|^{d/2}|a|\|h\|_{L^2(\Om_+)}.
\end{equation} 
Let $\{f_n\}$ be a sequence in $H^1_\eta(\Om_+)$ converging to $f$ in $H^1(\Om_+)$. 
We claim, and this is sufficient for our purposes, that $\EEeta f\in H^1(\Om)$ and  $\EEeta f_n\to \EEeta f$ in $H^1(\Om)$. Indeed, by assumption and \eqref{eq} it follows 
that $\EEeta f_n\to \EEeta f$ in $L^2(\Om)$. Also by assumption, for any $j=1,\ldots,d$, $\{\partial_j f_n\}$ is a Cauchy sequence in $L^2(\Om_+)$.
Hence,  by \eqref{eq} and \eqref{bbb}, $\{\partial_j( \EEeta f_n)\}$ is a Cauchy sequence in $L^2(\Om)$. Now, let  $\partial_j(\EEeta f_n)\to F_j$ in $L^2(\Om)$. 
It follows that  $\EEeta f\in H^1(\Om)$ and  $\partial_j (\EEeta f)=F_j$. The claim is proved.
\end{proof} 

It is worth mentioning here that obviously $C^\infty_c(\Omp)\subset H^1_\eta(\Omp)$ which implies that $H^1_\eta(\Omp)$ is dense in $L^2(\Omp)$. Moreover, the question when $H^1_\eta(\Om_+)$ 
are different for different $\eta$'s requires a comment. It occurs, that this heavily depends on the geometric relations of $\Om$ with the  hyperplanes $\langle\a_k\rangle^\bot$, $k=1,\ldots,m$. 
Namely, suppose that there exists $k\in\{1,\ldots,m\}$ such that  $B\setminus\langle\a_k\rangle^\bot\subset\Om$, where $B=B(x_0,r)$ is a Euclidean 
ball centered at $x_0\in\langle\a_k\rangle^\bot$. We can assume that $x_0\neq 0$ and $r$ is such that $B\cap \langle\a_j\rangle^\bot=\emptyset$ for $j\neq k$. Then, for 
$\eta,\eta'\in {\rm Hom} (W,\widehat{\mathbb Z}_2)$ such that $\eta(\s_{\a_k})\neq \eta'(\s_{\a_k})$, we have $H^1_\eta(\Omp)\neq H^1_{\eta'}(\Omp)$. Indeed, suppose that, for instance,  
$\eta(\s_{\a_k})=-1$ and  $\eta'(\s_{\a_k})=1$ and let $B$ be as above. Choose a $C^\infty$ function $\v$, symmetric in the hyperplane $\langle\a_k\rangle^\bot$, with support in $B$ and
equal to 1 on $\frac12 B$. It is easily seen that $\v^+\in H^1_{\eta'}(\Omp)$, while $\v^+\notin H^1_{\eta}(\Omp)$. 

Thus, for instance in the case when $\Om=\R$ and hence $\Omp=C_+$,  different $\eta$'s give different spaces $H^1_\eta(\Omp)$. This is a consequence of the fact that $\{\s_{\a_k}\}_1^m$ generate $W$.

Let $\mathfrak{t}=\mathfrak{t}_\Om$  denote the sesquilinear form defined on the Sobolev space $H^1(\Om)$,  
$$
\mathfrak{t}[F,G]=\int_{\Om}\sum_{j=1}^d \partial_jF(x)\,\overline{\partial_jG(x)}\,dx,\qquad F,G\in H^1(\Om).
$$
By $\mathfrak{t}_0$ we shall denote the  form $\mathfrak{t}$ restricted to $H^1_0(\Om)$. 
Analogously, let $\mathfrak{t}^+=\mathfrak{t}_{\Omp}$  denote the sesquilinear form defined on the Sobolev space $H^1(\Omp)$,
$$
\mathfrak{t}^+[f,g]=\int_{\Omp}\sum_{j=1}^d \partial_jf(x)\,\overline{\partial_jg(x)}\,dx, \qquad f,g\in H^1(\Omp), 
$$
and let $\mathfrak{t}^+_0$ denote the  form $\mathfrak{t}^+$ restricted to $H^1_0(\Omp)$. 
By $\mathfrak{t}^+_\eta$ we shall denote the sesquilinear form which is the restriction of $\mathfrak{t}^+$ to $H^1_\eta( \Omp)$,
$$
\mathfrak{t}^+_\eta[f,g]=\int_{\Omp}\sum_{j=1}^d \partial_jf(x)\,\overline{\partial_jg(x)}\,dx, \qquad f,g\in H^1_\eta(\Omp).
$$
Thus $\D(\mathfrak{t}^+_\eta)=H^1_\eta( \Omp)\subset H^1(\Omp)$.

Recall that in the general setting of a sesquilinear form $\mathfrak{s}$ with dense domain $\D(\mathfrak{s})$ and given on a Hilbert space $(\mathcal H, \langle\cdot,\cdot\rangle)$, 
the associated operator $A_\mathfrak{s}$ is defined by first determining its domain, 
$$
\D(A_\mathfrak{s})=\{h\in \D(\mathfrak{s})\colon \exists u_h\in \mathcal H\,\,\,\forall h'\in \D(\mathfrak{s})\,\,\, \mathfrak{s}[h,h']=\langle u_h,h'\rangle\},
$$
and then by setting its action on  $h\in \D(A_\mathfrak{s})$ by $A_\mathfrak{s}h=u_h$. 
If $\mathfrak{s}$ is Hermitian and closed, then $A_\mathfrak{s}$ is self-adjoint. 
If, in addition, $\mathfrak{s}$ is non-negative, then $A_\mathfrak{s}$ is non-negative.
See \cite[Chapter 10 and Section 3 of Chapter 12]{Sch}. Also recall, that closedness of a non-negative form $\mathfrak{s}$ means that the norm $\|x\|_{\mathfrak{s}}:=(\mathfrak{s}[x,x]^2+
\langle x,x \rangle^2)^{1/2}$ defined on $\D(\mathfrak{s})$ is complete.

Since the forms $\mathfrak{t}$, $\mathfrak{t}_0$, and $\mathfrak{t}^+$, $\mathfrak{t}^+_0$ and $\mathfrak{t}^+_\eta$, defined on $L^2(\Om)$ and $L^2(\Omp)$, respectively, are Hermitian, 
closed and non-negative, the associated operators are self-adjoint and non-negative. Clearly, closedness of $\mathfrak{t}$ and $\mathfrak{t}_0$,  and $\mathfrak{t}^+$ and $\mathfrak{t}^+_0$ 
is a consequence of the fact that the norms $\|\cdot\|_{H^1(\Om)}$ and $\|\cdot\|_{H^1(\Omp)}$ are complete and $H^1_0(\Om)$,  $H^1_0(\Omp)$ are closed subspaces in $H^1(\Om)$,  $H^1(\Omp)$, 
respectively. Closedness of  $\mathfrak{t}^+_\eta$ is a consequence of Proposition \ref{pro:clo}.

Let $-\Delta_{N,\Om}$ and $-\Delta_{D,\Om}$ be the Neumann and the Dirichlet Laplacians on $\Om$, the operators associated with $\mathfrak{t}$ and $\mathfrak{t}_0$, respectively.  
Let $-\Delta^+$, $-\Delta^+_0$ and $-\Delta^+_\eta$ denote the operators associated with $\mathfrak{t^+}$, $\mathfrak t^+_0$ and $\mathfrak{t^+_\eta}$, respectively.  
We shall call $-\Delta^+_\eta$ the  $\eta$-Laplacian on $\Omp$. It is clear that for $\eta\equiv1$ the corresponding operator is the Neumann  Laplacian on $\Omp$, while $-\Delta^+_0$ 
is the  Dirichlet  Laplacian on $\Omp$. Note, however, that if $\Om$ is such that $H^1(\Om)=H^1_0(\Om)$, then $-\Delta_{{\rm sgn}}^+$ coincides with $-\Delta^+_0$.

In addition, $-\Delta^+_0$ and each of the operators $-\Delta^+_\eta$ is indeed an extension of $-\Delta_{\Omp}$; we postpone explanation of this fact till the end of this section. 
Thus we have a number of natural self-adjoint extensions of the differential operator $-\Delta_{\Omp}$ with the Neumann Laplacian  and the Dirichlet Laplacian included.

The following result is crucial. It relates the domains $\D(-\Delta^+_\eta)$ and $\D(-\Delta_{N,\Om})$, or $\D(-\Delta^+_0)$ and $-\Delta_{D,\Om}$, as well as the actions of 
$-\Delta^+_\eta$ and $-\Delta_{N,\Om}$, or $-\Delta^+_0$ and  $-\Delta_{D,\Om}$, respectively. Recall that if  $\Om$ is such that $H^1(\Om)=H^1_0(\Om)$, then the case of Dirichlet Laplacians 
below, in Corollary \ref{cor:conv}, and in the proof of main theorem in Section \ref{sec:proof}, is absorbed by the general case.

\begin{proposition} \label{pro:main}
We have 
\begin{equation}\label{dom1}
\D(-\Delta^+_\eta)=\A_\eta (\D(-\Delta_{N,\Om}))^+
\end{equation}
and
\begin{equation}\label{for1}
(-\Delta^+_\eta)\big((\A_\eta F)^+\big)=\A_\eta(-\Delta_{N,\Om} F)^+, \qquad F\in \D(-\Delta_{N,\Om}).
\end{equation}
Analogous identities also hold when  $-\Delta^+_\eta$ and $-\Delta_{N,\Om}$ are replaced by $-\Delta_{D,\Omp}$ and $-\Delta_{D,\Om}$, respectively, 
and $\A_\eta $ is substituted by $\A_{\rm sgn}$. 
\end{proposition}
\begin{proof}
We first prove the inclusion $\subset$ in \eqref{dom1}. For brevity we write $-\Delta$ in place of $-\Delta_{N,\Om}$. Take $f\in\D(-\Delta^+_\eta)$. 
Hence $f\in H^1_\eta( \Om_+)$ that is $f=F^+$, where $F\in H^1(\Om)$ and $F=\A_\eta F$, and there is $u_f\in L^2(\Om_+)$ such that 
\begin{equation}\label{eqq1}
\forall g\in H^1_\eta( \Om_+)\quad\quad\mathfrak{t}^+_\eta[f,g]=\langle u_f,g\rangle_{L^2(\Om_+)},
\end{equation}
which also means that $(-\Delta^+_\eta)f=u_f$. Note that $\EEeta f=F$ and thus  $\EEeta f\in H^1(\Om)$. 
We shall verify that  it holds
\begin{equation}\label{eqq2}
\forall G\in H^1(\Om)\quad \quad\mathfrak{t}[\EEeta f,G]=\langle \EEeta(u_f),G\rangle_{L^2(\Om)},
\end{equation}
which means that  $\EEeta f\in \D(-\Delta)$ and hence  $f=(\EEeta f)^+=\A_\eta (\EEeta f)^+\in \A_\eta\big(\D(-\Delta^+_\eta)\big)^+$, 
and also that 
$$
(-\Delta)(\EEeta f)=\EEeta((-\Delta^+_\eta)f)\quad {\rm for}\quad f\in\D(-\Delta^+_\eta).
$$ 

For any $G\in H^1(\Om)$, using \eqref{bbb}, Lemma \ref{lem:first} and \eqref{H3} we have 
\begin{align*}
\mathfrak{t}[\EEeta f,G]=\sum_j\int_{\Om}\partial_j(\EEeta f)\,\overline{\partial_j G}
&=\sum_j\int_{\Om}\Big(\sum_i \mathcal{E}^{\eta\cdot g_{ji}}(\partial_i f)\Big)\,\overline{\partial_j G}\\
&=|W|\sum_i\int_{\Om_+} \partial_i f  \,\overline{\Big(\sum_j\A_{\eta \cdot g_{ji}}(\partial_j G)^+\Big)}\\
&=|W|\sum_i\int_{\Om_+} \partial_i f  \,\overline{\partial_i(\A_\eta G)^+}\\
&=|W|\mathfrak{t}^+_\eta[f,(\A_\eta G)^+].
\end{align*}
On the other hand,
$$
\langle \EEeta(u_f),G\rangle_{L^2(\Om)}=|W|\langle u_f,(\A_\eta G)^+\rangle_{L^2(\Om_+)},
$$
and hence inserting $(\A_\eta G)^+$ for $g$ in \eqref{eqq1} gives \eqref{eqq2}.

To prove the opposite inclusion in \eqref{dom1}, take $F\in\D(-\Delta)$. In particular, $F\in H^1(\Om)$ and, moreover, there is $U_F\in L^2(\Om)$ such that 
\begin{equation}\label{eq1}
\forall G\in H^1(\Om) \quad \quad \mathfrak{t}[F,G]=\langle U_F,G\rangle_{L^2(\Om)},
\end{equation}
which also means that $(-\Delta)F=U_F$. We shall verify that 
\begin{equation}\label{eq2}
\forall g\in H^1_\eta(\Om_+) \quad\quad   \mathfrak{t}^+_\eta[(\A_\eta F)^+,g]=\langle \A_\eta(U_F)^+,g\rangle_{L^2(\Om_+)},
\end{equation}
which means that $(\A_\eta F)^+\in \D(-\Delta^+_\eta)$ and completes the proof of \eqref{dom1}, and, in addition, 
that $(-\Delta^+_\eta)\big((\A_\eta F)^+\big)=\A_\eta((-\Delta)F)^+$ for $F\in\D(-\Delta)$, which completes the proof of \eqref{for1}.

Let $g\in H^1_\eta( \Om_+)$, where $g=G^+$ and $G=\A_\eta G$. Using \eqref{H3} and Lemma \ref{lem:first} we have 
\begin{align*}
\mathfrak{t}^+_\eta[(\A_\eta F)^+,g]&=\sum_j\int_{\Om_+}\partial_j((\A_{\eta} F)^+)\,\overline{\partial_j g}\\
&=\sum_j\int_{\Om_+}\Big(\sum_i\A_{\eta\cdot g_{ij}}(\partial_i F)^+\Big)\,\overline{\partial_j g}\\
&=|W|^{-d}\sum_i\int_{\Om} \partial_i F \,\overline{\Big(\sum_j\mathcal{E}^{\eta\cdot g_{ji}}(\partial_j g)\Big)}\\
&=|W|^{-d}\sum_i\int_{\Om} \partial_i F \,\overline{\partial_i(\EEeta g)}\\
&=|W|^{-d}\mathfrak{t}[F,\EEeta g].
\end{align*}
On the other hand,
$$
|W|\langle \A_\eta(U_F)^+,g\rangle_{L^2(\Om_+)}=\langle U_F,\EEeta g\rangle_{L^2(\Om)},
$$
and hence inserting $\EEeta g$ for $G$ in \eqref{eq1} gives \eqref{eq2}; note that $\EEeta g\in H^1(\Om)$ by assumption made on of $g$. 
This completes the proof of \eqref{eq2} and thus the conclusion following it, and hence finishes the proof of \eqref{dom1} and \eqref{for1}.  

We now comment the changes to be done in the above reasoning in order to prove the last claim of the proposition. We shall write $-\Delta$, $-\Delta_+$  
instead of $-\Delta_{D,\Om}$, $-\Delta_{D,\Om_+}$, respectively. For the inclusion $\subset$ in modified \eqref{dom1} we copy \textit{mutatis mutandis} 
the proof of $\subset$ in \eqref{dom1}. Clearly, we replace the symbol $H^1$ by $H^1_0$, note that $\EE^{{\rm sgn}} f\in H^1_0(\Om)$ provided 
$f\in H^1_0(\Omp)$, and substitute  $\EEeta$,  $\EE^{\eta\cdot g_{ji}}$,  $\A_\eta$,   $\A_{\eta\cdot g_{ji}}$, 
by $\EE^{{\rm sgn}}$, $\EE^{{\rm sgn}\cdot g_{ji}}$,  $\A_{\rm sgn}$,   $\A_{{\rm sgn}\cdot g_{ji}}$, respectively. The resulting identity is
$$
\mathfrak{t}[\EE^{{\rm sgn}} f,G]=|W|\mathfrak{t}^+[f,(\A_{\rm sgn} G)^+],
$$
which, together with
$$
\langle \EE^{{\rm sgn}}(u_f),G\rangle_{L^2(\Om)}=|W|\langle u_f,(\A_{\rm sgn} G)^+\rangle_{L^2(\Om_+)},
$$
shows an appropriate version of \eqref{eqq2} and its conclusion. Note that we used the fact that $(\A_{\rm sgn} G)^+\in H^1_0(\Omp)$. 
For the opposite inclusion in modified \eqref{dom1}, again 
we copy \textit{mutatis mutandis} the relevant part of the earlier reasoning with all previous symbol substitutions to conclude that 
$$
\mathfrak{t}^+[(\A_{\rm sgn} F)^+,g]=\frac1{|W|}\mathfrak{t}[F,\EE^{{\rm sgn}} g],
$$
and together with 
$$
|W|\langle \A_{\rm sgn}(U_F)^+,g\rangle_{L^2(\Om_+)}=\langle U_F,\EE^{{\rm sgn}} g\rangle_{L^2(\Om)},
$$
we end up with an appropriate version of \eqref{eq2} and the conclusions following it. 

This completes the proof of  modified \eqref{dom1} and \eqref{for1} and thus Proposition \ref{pro:main}.
\end{proof}

For further reference it is convenient to single out from the proof of Proposition \ref{pro:main} the following.
\begin{corollary} \label{cor:conv}
We have
\begin{equation*}
\EEeta\big({\rm Dom}(-\Delta_+^\eta)\big)\subset {\rm Dom}(-\Delta_{N,\Om})
\end{equation*}
and
\begin{equation*}
\EEeta\big((-\Delta_+^\eta)f\big)=(-\Delta_{N,\Om})(\EEeta f), \qquad f\in {\rm Dom}(-\Delta^+_\eta).
\end{equation*}
Analogous identities also hold when  $-\Delta^+_\eta$ and $-\Delta_{N,\Om}$ are replaced by $-\Delta_{D,\Omp}$ and $-\Delta_{D,\Om}$, respectively, 
and $\EEeta $ is substituted by $\mathcal E^{\rm sgn}$. 
\end{corollary}
We now come back to verification that $-\Delta_+^\eta$ and $-\Delta_0^+$ extend $-\Delta_{\Om_+}$. This will follow from definitions in term of forms, with an application of Gauss' formula 
for functions from Sobolev classes. Indeed, we claim that $C^2_c(\Om_+)\subset \D(-\Delta_+^\eta)$ and 
$(-\Delta_+^\eta)f=-\Delta_{\Om_+} f$, for $f\in C^2_c(\Om_+)$. To check this, take $f\in C^2_c(\Om_+)$ and let $\Omega_0$ be a bounded subset in $\Om_+$ with smooth boundary such that ${\rm supp}\,f\subset \Omega_0$. 
Recall that $\D(-\Delta_+^\eta)$ coincides with
$$
\{h\in \D(\mathfrak{t}^+_\eta)\colon \exists u_h\in L^2(\Om_+)\,\,\forall g\in \D(\mathfrak{t}^+_\eta)\,\, \mathfrak{t}^+_\eta[h,g]=\langle u_h,g\rangle_{L^2(\Om_+)}\},
$$
where $\mathfrak{t}^+_\eta[h,g]=\langle \nabla h,\overline{\nabla g}\rangle_{L^2(\Om_+)}$. Clearly $f\in\D(\mathfrak{t}^+_\eta)$, thus for every 
$g\in\D(\mathfrak{t}^+_\eta)$, since $f$ vanishes on $\partial\,\Omega_0$, by using  \cite[(D.4), Appendix D]{Sch} (note that $f,g\in H^1(\Om_+)$), we obtain
\begin{align*}
\mathfrak{t}^+_\eta[f,g]&=\int_{\Om_+}\Big(\sum_j\partial_j f\,\overline{\partial_j g}\Big)\,dx
=\int_{\Omega_0}\Big(\sum_j\,\partial_j f\overline{\partial_j g}\Big)\,dx\\
&=\int_{\Omega_0}(-\Delta_{\Om_+})\,f\overline{g}\,dx
=\langle (-\Delta_{\Om_+})f,g\rangle_{L^2(\Om_+)}.
\end{align*}
This proves the claim. 

Finally, the question when $-\Delta_\eta^+$ are different for different $\eta$'s also requires a comment. We focus on looking at the domains 
${\rm Dom}(-\Delta_+^\eta)$ and apply an argument similar to that used in the analysis of differences between the spaces $H^1_\eta(\Omp)$. 
As in that case the geometric relations of $\Om$ with the  hyperplanes $\langle\a_k\rangle^\bot$, $k=1,\ldots,m$, are essential. Therefore, 
suppose that there exists $k\in\{1,\ldots,m\}$ such that  $B\subset\Om$, where $B$ is a Euclidean ball centered at $0\neq x_0\in\langle\a_k\rangle^\bot$ 
and $B\cap \langle\a_j\rangle^\bot=\emptyset$ for $j\neq k$. Then for $\eta,\eta'\in {\rm Hom} (W,\widehat{\mathbb Z}_2)$ such that 
$\eta(\s_{\a_k})\neq \eta'(\s_{\a_k})$, we have ${\rm Dom}(-\Delta_+^\eta)\neq {\rm Dom}(-\Delta_+^{\eta'})$. Indeed, suppose that, for instance,  
$\eta(\s_{\a_k})=-1$ and  $\eta'(\s_{\a_k})=1$ and let $B$ be as above. Choose a $C^\infty$ function $\v$  with support in $B$ and
equal to 1 on $\frac12 B$, which is $\eta'$-symmetric. Then $\v\in {\rm Dom}(-\Delta_{N,\Om})$ and, consequently, by \eqref{dom1}, 
$\v^+\in {\rm Dom}(-\Delta_+^{\eta'})$. On the other hand it is clear that  $\v^+\notin H^1_{\eta}(\Omp)$ and hence $\v^+\notin {\rm Dom}(-\Delta_+^\eta)$. 

Thus, for instance in the case of $\Om=\R$ and hence $\Omp=C_+$,  different $\eta$'s give different operators $-\Delta_+^\eta$. Again this is a 
consequence of the fact that $\{\s_{\a_k}\}_1^m$ generate $W$.

\section{Proof of the main result} \label{sec:proof}
The following \textit{commuting property} of the spectral functional calculus is well known: if $A$ is a self-adjoint operator on a Hilbert space 
$\mathcal H$ and $B\colon \mathcal H\to \mathcal H$ is a bounded operator such that $BA\subset AB$, then also $B\Psi(A)\subset \Psi(A)B$, for any Borel 
function $\Psi$ on $\mathbb R$. At least as a folklore,  the following two-Hilbert space and two-operator version of this commuting property  is also known: 
if $A_1$ and $A_2$ are self-adjoint operators on  Hilbert spaces $\mathcal H_1$ and $\mathcal H_2$, respectively, and $B\colon \mathcal H_1\to \mathcal H_2$ 
is a bounded operator such that $BA_1\subset A_2B$, then also $B\Psi(A_1)\subset \Psi(A_2)B$, for any Borel function $\Psi$ on $\mathbb R$; see \cite{MS} for additional comments. We shall refer to this property as to the \textit{intertwining property}. Recall, that the latter inclusion precisely means the inclusion 
of domains (which is equivalent to the statement that $B({\rm Dom}(\Psi(A_1)))\subset {\rm Dom}(\Psi(A_2))$) and the identity 
$B\big(\Psi(A_1)x\big)=\Psi(A_2)\big(Bx\big)$ for all $x\in{\rm Dom}(\Psi(A_1))$.

Passing to the proof of Theorem \ref{thm:main} consider first the case of the $\eta$-Laplacians on $\Omp$ and prove \eqref{kerN}. As in the proof of 
Proposition \ref{pro:main}, for brevity, if not otherwise stated, we write $-\Delta$ instead of $-\Delta_{N,\Om}$,  and consequently, $\Psi(-\Delta)$  
instead of $\Psi(-\Delta_{N,\Om})$. Analogously, we write $K^\Psi$  rather than $K^\Psi_{-\Delta_{N,\Om}}$. 

It is intuitively clear and easily seen by using the definition of $-\Delta$ in terms of the sesquilinear form $\mathfrak{t}$, that for any $g \in W$, 
the operator $T_g \colon F\mapsto F_g$, $F_g (x)=F(g  x)$, which is bounded on $L^2(\Om)$, commutes with $-\Delta$. This implies that 
$$
\big(\Psi(-\Delta)F\big)_g =\Psi(-\Delta)(F_g ), \qquad F\in \D\,\big(\Psi(-\Delta)\big).
$$
Consequently, since $\D(\Psi(-\Delta))$ is dense in $L^2(\Om)$, the kernel $K^\Psi$  satisfies
\begin{align}\label{symm}
K^\Psi(g  x,y)=K^\Psi(x,g^{-1}y),\qquad (x,y)\in\Om\times\Om\,-\, a.e.
\end{align}

On the other hand, \eqref{for1} of Proposition \ref{pro:main} says that the operator $\A^+_\eta\colon F\mapsto (\A_\eta F)^+$, which is bounded from $L^2(\Om)$ 
into $L^2(\Omp)$, intertwines the operators $-\Delta$ and $-\Delta^+_\eta$ in the sense that $\A^+_\eta\circ (-\Delta)\subset (-\Delta^+_\eta)\circ \A^+_\eta$. 
It follows, by applying  the beforementioned  intertwining property to $\mathcal H_1=L^2(\Om)$ and $\mathcal H_2=L^2(\Omp)$, $A_1=-\Delta$ and 
$A_2=-\Delta_\eta^{+}$, and $B=\A^+_\eta$, that also $\A^+_\eta$ intertwines the operators $\Psi(-\Delta)$ and $\Psi(-\Delta^+_\eta)$; in particular, 
\begin{equation}
\label{delta}\A_\eta\big(\Psi(-\Delta)F\big)^+=\Psi(-\Delta^+_\eta)\big((\A_\eta F)^+ \big), \quad {\rm for}\quad F\in \D\big(\Psi(-\Delta)\big).
\end{equation}

On the other side, Corollary \ref{cor:conv} says that the operator $\EEeta \colon f\mapsto \EEeta  f$, which is bounded from $L^2(\Omp)$ into $L^2(\Om)$, 
intertwines the operators $-\Delta^+_\eta$ and $-\Delta$. It follows that  $\EEeta $ also intertwines the operators $\Psi(-\Delta^+_\eta)$ and 
$\Psi(-\Delta)$ and hence, 
\begin{equation}\label{delta2}
\EEeta \big( {\rm Dom}(\Psi(-\Delta^+_\eta))\big)\subset  \D\big(\Psi(-\Delta)\big).
\end{equation}
Here the intertwining property was applied to $\mathcal H_1=L^2(\Omp)$ and  $\mathcal H_2=L^2(\Om)$, $A_1=-\Delta^+_\eta$ and $A_2=-\Delta$, and $B=\EEeta $.

Now take $f\in \D(\Psi(-\Delta^+_\eta))$ and let $F=\EEeta  f$. Then, by \eqref{delta2},  $F\in \D(\Psi(-\Delta))$ and, since $F$ is $\eta$-invariant, 
$f=F^+=(\A_\eta F)^+$. Hence, by \eqref{delta},  for $x\in \Omp$, we have
\begin{align*}
\Psi(-\Delta^+_\eta)f(x)
&=\frac1{|W|}\sum_{g \in W}\eta(g)\Psi(-\Delta)F(gx)\\
&=\frac1{|W|}\sum_{g \in W}\eta(g)\int_{\Om} K^\Psi(gx,y)F(y)\,dy\\
&=\frac1{|W|}\sum_{g \in W}\eta(g)\sum_{g'\in  W}\eta(g')\int_{\Omp} K^\Psi(gx,g'y)f(y)\,dy,\\
&=\int_{\Omp}\Big(\sum_{g''\in W} \eta(g'') K^\Psi(g''x,y)\Big)f(y)\,dy,
\end{align*}
where, for the last identity, we used \eqref{symm} to obtain
\begin{align*}
\sum_{g \in W}\eta(g)\sum_{g '\in W}\eta((g')^{-1}) K^\Psi(g  x,g'y)&=\sum_{g \in W}\eta(g)\sum_{g'\in  W}\eta((g')^{-1}) K^\Psi((g')^{-1}g  x,y)\\
&=|W|\sum_{g''\in W }\eta(g'') K^\Psi(g''x,y).
\end{align*}
This means that $\Psi(-\Delta^+_\eta)$ has an integral kernel and \eqref{kerN} holds. 

We now comment the changes to be done in the above reasoning in order to prove  \eqref{kerD}. For brevity, we now write $-\Delta$, $-\Delta_+$, 
$\Psi(-\Delta)$, $\Psi(-\Delta_+)$ and  $K^\Psi$ instead of $-\Delta_{D,\Om}$, $-\Delta_{D,\Om_+}$, $\Psi(-\Delta_{D,\Om})$, $\Psi(-\Delta_{D,\Om_+})$, 
$K^\Psi_{-\Delta_{D,\Om}}$, respectively. Then we repeat the arguments we used in the previous case so that we again obtain \eqref{symm} and by using 
version of \eqref{for1} of Proposition \ref{pro:main} with $\eta={\rm sgn}$ (and $-\Delta$ replacing $-\Delta_{D,\Om_+}$) we arrive to a version of 
\eqref{delta} with $\A_\eta$ replaced by $\A_{\rm sgn}$. Consequently, picking up $f\in \D(\Psi(-\Delta_+))$ and  $F\in \D(\Psi(-\Delta))$ such that 
$(\A_{\rm sgn} F)^+=f$ and considering $F=\EE^{\rm sgn} f$, we follow the main calculation of the previous case and end up with
$$
\Psi(-\Delta_+)f(x)=\int_{\Om_+}\Big(\sum_{g''\in W} {\rm sgn}(g'')K^\Psi(g''x,y)\Big)f(y)\,dy, \qquad x\in\Omp.
$$
The proof of Theorem \ref{thm:main} is completed.

\section{Examples} \label{sec:app}
\subsection{Groups associated with orthogonal root systems} \label{sec:ortho}
This example is related to the one discussed in \cite[Section 4.1]{MS}. There, for an orthogonal root system and an arbitrary open and approprietely symmetric $\Omega\subset\R$, an inductive argument was applied to obtain required relations between relevant kernels. Here we focus on  $\Omega=\R$ and consider a Weyl chamber coming from an orthogonal root system. Then we apply directly our main theorem to obtain relations between kernels labeled by a variety of relevant homomorphisms. In addition, properties of $\eta$-heat kernels are discussed together with their application to solving some mixed Neumann-Dirichlet initial-boundary value problems.

Let $R$ be  an orthonormal root system in $\R$ with a chosen subsystem $R_+$  of positive roots (orthonormality of $R$ means that  $R_+$ is orthonormal 
as a set of vectors). Without any loss of generality, possibly by rotating and reflecting the coordinate axes, we can assume that $R_+=\{e_{j_1},\ldots,e_{j_k}\}$, where $1=j_1<j_2<\ldots<j_k\le d$, and $e_j$ is the $j$th coordinate unit vector. Thus, given $1\le k\le d$ and $J=(j_1,\ldots,j_k)$ as above, let 
$R_+^{(k,J)}=\{e_{j_s}\colon s=1,\ldots,k\}$ be the system of positive roots so that $R^{(k,J)}=R_+^{(k,J)}\sqcup(-R_+^{(k,J)})$ is the orthonormal root system 
in $\R$. In what follows to fix the attention we consider only the case when $k=d$; the other cases, $1\le k\le d-1$, can be treated analogously. 

The positive Weyl chamber then is $\R_+=(0,\infty)^d$ and we identify the Weyl group corresponding to $R^{(d,J_d)}$, $J_d=(1,2,\dots,d)$, with 
$\widehat{\mathbb Z}_2^d$, and the action of $\ve\in \widehat{\mathbb Z}_2^d$ on $\R$ is through $x\to \ve x=(\ve_1x_1,\ldots,\ve_dx_d)$. Next, we identify 
${\rm Hom} (\widehat{\mathbb Z}_2^d, \widehat{\mathbb Z}_2)$ with $\mathbb Z^d_2$, where, this time, $\mathbb Z_2=\{0,1\}$ with addition modulo 2. In this identification a homomorphism $\eta \in \mathbb Z^d_2$ of the Weyl group represented by $\widehat{\mathbb Z}_2^d$ into $\widehat{\mathbb Z}_2$ acts  through 
$\ve \to \ve^\eta:=\prod_{j=1}^d \ve_j^{\eta_j}$ for  $\ve=(\ve_1,\ldots,\ve_d)$. 

Also, in these identifications, the $\eta$-symmetrization operator is 
$$
\A_\eta F(x)=\frac1{2^d}\sum_{\ve\in\{-1,1\}^d} \ve^\eta F(\ve x), \qquad x\in\R.
$$
Moreover, \eqref{H3} simplifies to the form $\partial_j(\A_\eta F)=\A_{\tau_j(\eta)}(\partial_j F)$, where for $\eta=(\eta_1,\ldots,\eta_d)\in\mathbb Z_2^d$ 
and $j\in\{1,\ldots,d\}$ we write $\tau_j(\eta)=(\eta_1,\ldots,1-\eta_j,\ldots,\eta_d)$ (this means the '$\mathbb{Z}_2$\,-reflection' on the $j$th coordinate). 

Let  $\{\exp(-t(-\Delta))\}_{t>0}$ denote  the \textit{heat semigroup} associated with $-\Delta$, the Neumann Laplacian on $\R$ (that coincides with the 
Dirichlet Laplacian). For every $t>0$, $\exp(-t(-\Delta))$ is a bounded on $L^2(\R)$ operator with kernel $p_t^{-\Delta}(x,y)=p_t^{(d)}(x-y)$, $x,y\in\R$, where
$$
p_t^{(d)}(w)=(4\pi t)^{-d/2}\exp\big(-\|w\|^2/4t\big), \qquad w\in\R,
$$
is the $d$-dimensional Gauss-Weierstrass kernel.

Similarly, let $\{\exp(-t(-\Delta^+_\eta))\}_{t>0}$ denote  the $\eta$-\textit{heat semigroup} associated with the $\eta$-Laplacian $-\Delta^+_\eta$, 
the semigroup of bounded operators on $L^2(\Rp)$ and let $p_t^{\eta,+}$ denote the corresponding kernel. Due to the property 
$p_t^{(d)}(w)=\prod_{j=1}^d p_t^{(1)}(w_j)$, as an immediate consequence of Corollary \ref{cor:main} we obtain the following.
\begin{corollary} \label{cor:cor1}
For every $t>0$, $\exp\big(-t(-\Delta^+_\eta)\big)$ is an integral operator in $L^2(\Rp)$ with kernel given by
\begin{equation*}
p_t^{\eta,+}(x,y)=\prod_{j=1}^d\big(p_t^{(1)}(x_j-y_j)+(-1)^{\eta_j} p_t^{(1)}(x_j+y_j) \big), \qquad x,y\in \Rp.
\end{equation*}
\end{corollary}

Note that (expected) properties of the $\eta$-\textit{heat kernel} $\{p_t^{\eta,+}\}_{t>0}$, analogous to these that correspond to the Neumann and 
Dirichlet heat kernels on open subsets of $\R$, readily follow by inspection. Indeed, smoothness of $p_t^{\eta,+}(x,y)$ jointly in the variables $t>0$, 
$x,y\in \Rp$, is obvious. Equally obvious is symmetry with respect to the $x,y\in \Rp$ variables. Also the fact that for every fixed $y\in \Rp$ the function 
$v(t,x):=p_t^{\eta,+}(x,y)$ satisfies the equation $(\partial_t-\Delta_{\Rp})v=0$ on $(0,\infty)\times \Rp$ as well as strict positivity of $p_t^{\eta,+}(x,y)$ 
for $x,y\in \Rp$, is easily seen. For the qualitative properties of $p_t^{\eta,+}(x,y)$ one easily observes that for every fixed $x\in \Rp$ and $t>0$ it holds:

\noindent a)  
$$
0<\int_{\Rp}p_t^{\eta,+}(x,y)\,dy\le1, \qquad \lim_{t\to0^+}\int_{\Rp}p_t^{\eta,+}(x,y)\,dy=1;
$$
moreover, for $\eta=\textbf{0}$ the estimated integral equals 1; otherwise for every fixed $x\in \Rp$ and $t>0$ the integral is strictly less than 1.

\noindent b) for every $\epsilon>0$ and $0<\delta<{\rm dist}(x, \partial \Rp)$ there exist $t_0>0$ such that for all $0<t<t_0$ one has
$$
\int_{\Rp\setminus B(x,\delta)}p_t^{\eta,+}(x,y)\,dy\le\epsilon,
$$
where $B(x,\delta)$ denotes the Euclidean ball centered at $x$ with radius $\delta$.

Let us focus, therefore, on the "boundary" properties. First, notice that the boundary of $\Rp$ consists of $d$ facets 
$\mathcal{F}_j=\{x\in \R\colon x_j=0\,\,\,{\rm and}\,\,\, x_k\ge0\,\,{\rm for}\,\, k\neq j\}$, $j=1,\ldots,d$. It is now easy to check that $p_t^{\eta,+}(x,y)$ satisfies mixed Neumann-Dirichlet conditions:  for every fixed $y\in \Rp$, the normal outward derivative of $p_t^{\eta}(\cdot,y)$ vanishes on any facet 
$\mathcal{F}_j$ (to be precise on the 'interior' of $\mathcal{F}_j$, where the derivative exists), where $\eta_j=0$, and $p_t^{\eta}(\cdot,y)$ vanishes 
on any facet $\mathcal{F}_j$, where $\eta_j=1$. Clearly, this agrees with the Neumann boundary conditions when $\eta=\textbf{0}$ or  with the 
Dirichlet boundary conditions when $\eta=\textbf{1}$.

To simplify the notation, in what follows, given $\eta\in \mathbb Z^d_2$ we split the boundary $\partial \Rp$ into $(\partial \Rp)_{\eta,N}$ and 
$(\partial \Rp)_{\eta,D}$, where the former part includes all facets $\mathcal{F}_j$, where $\eta_j=0$, and the latter part includes all facets 
$\mathcal{F}_j$, where $\eta_j=1$. We now come to the following initial-boundary problem.
$$
(\partial_t-\Delta_{\Rp})u(t,x)=0, \qquad (t,x)\in(0,\infty)\times \Rp,
$$
$$
u(0,x)=f(x), \qquad x\in \Rp,
$$
\begin{equation}\label{eta}
\frac{\partial u}{\partial \nu}(t,\overline{x})=0\,\,{\rm for}\,\,\overline{x}\in (\partial \Rp)_{\eta,N},\qquad u(t,\overline{x})=0\,\,{\rm for}
\,\,\overline{x}\in (\partial \Rp)_{\eta,D},\quad t>0.
\end{equation}
Here $f$ is a given (suitable) function on $\Rp$ and $\nu$ denotes the normal outward vector on $\partial \Rp$ (whenever exists). 
As one can expect the solution to this problem is given by
\begin{equation}\label{sol}
u(t,x)=\int_{\Rp}  p_t^{\eta,+}(x,y)f(y)\,dy, \qquad (t,x)\in(0,\infty)\times \overline{\Rp}.
\end{equation}
Note that $p_t^{\eta,+}(\overline{x},y)$, $y\in \Rp$, and thus also $u(t,\overline{x})$ are well defined for $\overline{x}\in \partial \Rp$.

A typical result that can be reached by employing clasical arguments and using beforementioned properties of the $\eta$-heat kernel reads as follows 
(we leave  details of the proof for interested readers).
\begin{proposition} \label{pro:inibound}
Assume that $f$ is continuous and bounded on $\overline{\Rp}$ and let $u$ be given by \eqref{sol}. Then

\noindent 1) $u$ is $C^\infty$ on $(0,\infty)\times \Rp$;

\noindent 2) $u$ satisfies the equation $(\partial_t-\Delta_{\Rp})u(t,x)=0$ on $(0,\infty)\times \Rp$;

\noindent 3) $u$ satisfies the initial condition $u(0,x)=f(x)$, $x\in \Rp$, in the sense that
$$\lim_{t\to0^+}u(t,x)=f(x)\quad for\,\, every \quad x\in \Rp;$$

\noindent 3) $u$ satisfies the $\eta$-boundary conditions \eqref{eta}.
\end{proposition}

\subsection{Dihedral groups} \label{sec:exs}
Special cases (dyadic cones) of the examples considered in this subsection were discussed in \cite[Section 4.4.2]{MS}.

Let $D_n$, $n\ge3$, be the reflection  group associated with the root system $R=\{z_j\colon j=0,\ldots,2n-1\}$ in $\RR\simeq \C$, where $z_j=e^{i\pi j/n}$. 
$D_n$ is called the dihedral group and geometrically this is the group of isometries of the regular $n$-gon centered at the origin with one of the vertices 
located at $i\simeq(0,1)\in\RR$. $D_n\subset O(2)$ is isomorphic with the semidirect product $\mathbb Z_n \rtimes \widehat{\mathbb Z}_2$ and $|D_n|=2n$. 
In this isomorphism $(1,1)$ corresponds to $r$, the counterclockwise rotation by $2\pi/n$ around the origin, while $(0,-1)$ corresponds to $\s$, the 
reflection with respect to the line orthogonal to $z_0$. Moreover,
$$
D_n=\{1,r,\ldots,r^{n-1},\s,r\s,\ldots,r^{n-1}\s\}.
$$
and $1,r,\ldots,r^{n-1}$ are all rotations in $D_n$, while $\s,r\s,\ldots,r^{n-1}\s$ are all reflections in $D_n$. 

The structure of normal subgroups of index 2 in $D_n$ is known. Namely,  $\NN_0=\langle r\rangle$ is the unique such a normal subgroup for $n$ odd, 
while for $n$ even, apart of $\NN_0$, there are two additional, $\NN_1=\langle r^2,\s\rangle$ and $\NN_2=\langle r^2,r\s\rangle$. Clearly, if $\NN$ 
is any normal subgroup of index 2 in $D_n$, then
$$
\eta_\NN(g)= \left\{ \begin{array}{ll}
1 & \textrm{for $g\in \NN$},\\
-1 & \textrm{for $g\notin \NN$,}
\end{array} \right.
$$
is a homomorphism of $D_n$ onto $\widehat{\mathbb Z}_2$ with ${\rm ker}(\eta_\NN)=\NN$, and every nontrivial homomorphism of $D_n$ into $\widehat{\mathbb Z}_2$ 
is of this form. In this notation, $\eta_{\NN_0}= {\rm sgn}$. Therefore, the complete list of homomorphisms of $D_n$ into $\widehat{\mathbb Z}_2$ is as follows.
For $n$ odd the only such homomorphisms are $\textbf{1}$ and ${\rm sgn}$, while for $n$ even there are two additional, $\eta_{\NN_1}$ and $\eta_{\NN_2}$.

We shall write down the explicit form of the Neumann and Dirichlet heat kernels, as well as the two exceptional $\eta$-heat kernels,  
for the open cone on the plane with vertex at the origin and aperture $\pi/n$. This cone,
$$
C_+^{(n)}=
\begin{cases}
\{\rho e^{i\theta}\colon \rho>0,\,\,\,\,0<\theta<\pi/n\},& n\,\,\,\, {\rm even},\\
\{\rho e^{i\theta}\colon \rho>0,\,\,\,\,|\theta|<\pi/2n\},& n\,\,\,\, {\rm odd},
\end{cases}
$$
is the positive Weyl chamber for the following selection of the set of positive roots:
if $n=2k$, $k\ge2$, then $R_+=\{z_{3k+1},z_{3k+2},\ldots, z_{4k-1},z_0,z_1,\ldots, z_{k}\}$ with $\Pi=\{z_{3k+1},z_{k}\}$ as the set of simple roots; 
if $n=2k+1$, $k\ge1$, then $R_+=\{z_{3k+2},\ldots, z_{4k+1},z_0,\ldots, z_{k}\}$ with $\Pi=\{z_{3k+2},z_{k}\}$ as the set of simple roots.

 For brevity we now write $p_t$ for the two dimensional Gauss-Weierstrass kernel $p_t^{(2)}$.

The Neumann heat kernel on $C_+^{(n)}$, which is the $\eta$-heat kernel for $\eta\equiv 1$,  is given by
$$
p_t^{N,\,\pi/n}(x,y)=\sum_{g\in D_n} p_t(gx-y)=\sum_{m=0}^{n-1}\big(p_t(r^mx-y)+p_t(r^m\s x-y)\big);
$$
here $x,y\in C_+^{(n)}$. The Dirichlet heat kernel on $C_+^{(n)}$,  which is the $\eta$-heat kernel for $\eta={\rm sgn}$ (and corresponds to 
${\rm ker}\,\eta=\NN_0$), is for $x,y\in C_+^{(n)}$ given by
$$
p_t^{D,\,\pi/n}(x,y)=\sum_{g\in D_n} {\rm det}(g)\,p_t(gx-y)=\sum_{m=0}^{n-1}\big(p_t(r^mx-y)-p_t(r^m\s x-y)\big).
$$

It is worth mentioning here that the formula for $p_t^{D,\,\pi/n}(x,y)$ in terms of a convergent series is well known (in fact any aperture is admitted) 
and goes back to an old formula of Carslaw and Jaeger. See \cite{BS}, where this formula was generalized to higher dimensions (and \cite{MS}, where the 
analogous formula is stated in the Neumann case). The formula we have in mind, for $n$ even reads: in polar coordinates, for 
$x=\rho e^{i\theta}\in C_+^{(n)}$, $y=re^{i\xi}\in C_+^{(n)}$,
$$
p_t^{D,\,\pi/n}(x,y)=\frac n{2\pi t}\exp\left(-\frac{\rho^2+r^2}{4t}\right)\sum_{j=1}^\infty I_{jn}\Big(\frac{\rho r}{2t}\Big) 2\sin(jn\theta)\sin(jn\xi).
$$
Here $I_\nu(z)$ denotes the modified Bessel function of order $\nu$ (in the Neumann case  sines are replaced by cosines). For $n$ odd $\theta$ and $\xi$ have to be replaced by $\theta+\pi/2n$ and $\xi+\pi/2n$; for such $n$, $C_+^{(n)}$ rotated counterclockwise by $\pi/2n$ sticks to the form of $C_+^{(n)}$ for $n$ even.

In the exceptional cases, when $n=2k$, $k\ge2$, we have
$$
\NN_1=\{1,r^2,\ldots,r^{2(k-1)},\s,r^2\s,\ldots,r^{2(k-1)}\s\},
$$
and 
$$
\NN_2=\{1,r^2,\ldots,r^{2(k-1)},r\s,r^3\s,\ldots,r^{2k-1}\s\}.
$$
Therefore the corresponding $\eta_j$-heat kernels, $j=1,2$, $x,y\in C_+^{(n)}$, are
\begin{align*}
&\,\,\,\,\,\,\,\,\,p_t^{\eta_1,\,\pi/n}(x,y)\\
&=\sum_{g\in \NN_1}\big(p_t(gx-y)-p_t(rg x-y)\big)\\
&=\sum_{m=0}^{k-1}\big(p_t(r^{2m}x-y)-p_t(r^{2m+1} x-y)+p_t(r^{2m}\s x-y)-p_t(r^{2m+1}\s x-y)\big), 
\end{align*}
and
\begin{align*}
&\,\,\,\,\,\,\,\,\,p_t^{\eta_2,\,\pi/n}(x,y)\\
&=\sum_{g\in \NN_2}\big(p_t(gx-y)-p_t(g\s x-y)\big)\\
&=\sum_{m=0}^{k-1}\big(p_t(r^{2m}x-y)-p_t(r^{2m+1} x-y)+p_t(r^{2m+1}\s x-y)-p_t(r^{2m}\s x-y)\big).
\end{align*}

\subsection*{Acknowledgements} Research was supported by funds of Faculty of Pure and Applied Mathematics, Wroc\l{}aw University of Science and Technology, 
\# 0401/0155/18. The author thanks Jacek Ma\l{}ecki for  discussions.

\end{document}